\documentclass{article}

\usepackage{graphicx} % more modern
\usepackage{caption}
\usepackage{subcaption}
\usepackage{epstopdf}
\usepackage{amssymb}
\usepackage{algorithm}
\usepackage{algpseudocode}
\usepackage{amsmath}
\usepackage{amsthm}
\usepackage{amsfonts}
\usepackage{mathtools}
\usepackage{color}
\usepackage{indentfirst}
\usepackage{bm}
\usepackage[margin=3.0cm]{geometry}
\usepackage{authblk}

\renewcommand{\i}{\pmb{\mathrm{i}}}
\newtheorem{theorem}{Theorem}
\newtheorem{lemma}{Lemma}
\newtheorem{Remark}{Remark}

\graphicspath{{./Figures/}}

\title{The discrete sign problem: uniqueness, recovery algorithms and phase retrieval applications}
\author[1]{Ben Leshem}
\author[2]{Oren Raz}
\author[3]{Ariel Jaffe}
\author[3]{Boaz Nadler}

\affil[1]{Department of Physics of Complex Systems, Weizmann Institute of Science, Rehovot 76100, Israel}
\affil[2]{Department of Chemistry and Biochemistry, University of Maryland, College Park, MD 20742, U.S.A.}
\affil[3]{Department of Computer Science 
and Applied Mathematics , Weizmann Institute of Science, Rehovot 76100, Israel}
\date{}
\begin{document}

\maketitle
\begin{abstract}
In this paper we consider the following real-valued and finite dimensional specific instance of the 1-D classical phase retrieval problem. Let ${\bf F}\in\mathbb{R}^N$ be an $N$-dimensional vector, whose discrete Fourier transform has a compact support. The sign problem is to recover ${\bf F}$ from its magnitude $|{\bf F}|$.\ First, in contrast to the classical 1-D phase problem which in general has multiple solutions, we prove that with sufficient over-sampling, the sign problem admits a unique solution. Next, we show that the sign problem can be viewed as a special case of a more general piecewise constant phase problem. Relying on this result,  we derive a computationally efficient and robust to noise sign recovery algorithm. In the noise-free case and with a sufficiently high sampling rate,  our algorithm is guaranteed to recover the true sign pattern. Finally, we present two phase retrieval applications of the sign problem: (i) vectorial phase retrieval with three measurement vectors; and (ii) recovery of two well separated 1-D objects.    
\end{abstract}


\section{Introduction}\label{Section:Intro}

The recovery of a signal from modulus (absolute value) measurements of its Fourier
transform, known as {\em phase retrieval}, is a classical problem with a broad range of applications, including X-ray crystallography \cite{Millane_90}, astrophysics \cite{fienup1987phase}, lensless imaging \cite{miao1998phase,Lensless_NatPhys_2006}, and characterization of ultra-short pulses \cite{Walmsley_Review}, to name but a few.  

From a mathematical perspective, fundamental questions regarding the uniqueness
of the phase problem and development of reconstruction algorithms have  been topics of intense research for several decades. Extensive works have studied the number of possible solutions to the phase problem both in one-dimensional and two-dimensional settings, typically under the assumption that the underlying signal has a compact support,  see for example 
\cite{balan2015reconstruction,Sodin,burge1976phase,Klibanov_95,Sanz} and many additional references therein. For a recent survey and new results on the uniqueness of the 1-D phase problem, see Beinert and Plonka \cite{beinert2015ambiguities}. From the computational aspect, as the classical phase problem is non-convex, the most commonly used phase-retrieval methods are iterative \cite{Bauschke_02,elser2003phase,Fienup,Review_Marchesini}.  These may require careful initialization, may exhibit limited robustness to noise in various settings,  and in general, even in the absence of noise, are not guaranteed to converge to the correct solution.

In recent years there has been a renewed surge of interest in the phase problem. This was motivated in part by proposals of new measurement schemes coupled with novel convex-optimization approaches that  lead to strong guarantees on correct recovery. Examples include coded diffraction patterns, polarization type schemes and other methods involving multiple illuminations \cite{candes2015phase,Balan_09,alexeev2014phase,OurIEEE,novikov2015illumination},  matrix completion approaches to phase retrieval via semi-definite programs \cite{waldspurger2015phase,Candes} and sparsity-based recovery methods \cite{sidorenko2015sparsity,Eldar_Segev_Nature}.

Motivated by several phase retrieval applications, in this paper we consider a finite dimensional and real-valued particular instance of the general 1-D phase problem, which we denote as the {\em sign problem}. With precise definitions appearing in Section \ref{Section:Setup}, its formulation is as follows: Let \({\bf F}\) be an $N$-dimensional real-valued vector, whose discrete Fourier transform, ${\bf f}= \mathcal{DFT}\{{\bf F}\}\in\mathbb{C}^N$,
has a compact support of length $\tau+1$. The sign problem is to recover the sign pattern \({\bf s}=sign({\bf F})\in\{\pm 1\}^N\) from possibly noisy measurements of  its magnitude $|{\bf F}|$.

In this paper we perform a detailed study of the finite dimensional sign problem, including its uniqueness and the development of a stable reconstruction algorithm. We also present its application to two practical phase problems. Specifically we make the following contributions. First, in Section \ref{Section:MathUniq}, Theorem \ref{Theorem:SignUniq} we prove that if \(N>2\tau\), our discrete sign problem admits a unique solution, up to a global $\pm 1$ sign ambiguity. Our proof technique, similar to \cite{OurIEEE,Sodin,beinert2015ambiguities}, is based on analyzing the roots of high degree polynomials. 
Since finding the roots of high degree polynomials is known to be an ill-conditioned problem \cite{trefethen1997numerical}, our above proof does not directly lead to a stable reconstruction algorithm. 

To robustly solve the sign problem we thus consider a different approach:\ First, we study the structure of the solutions to the sign problem. In Lemma \ref{lem:sign_changes}, we show that the sign pattern ${\bf s}=sign({\bf F})\in\{\pm 1\}^N$ cannot be arbitrary, but rather has at most $\tau$ sign changes. Next, we relax the constraint that ${\bf s}\in\{\pm 1\}^N$ and allow the sign pattern to be a complex-valued \(N\) dimensional vector, which guided by Lemma \ref{lem:sign_changes}, is piecewise constant over at most \(\tau+1\) intervals. This leads us to study the following two questions: (i) is it  possible to detect at least parts of these intervals, where the underlying sign pattern is constant? and (ii) does such an {\em over-segmentation} of the indices $1,\ldots,N$ to intervals of constant values indeed retains uniqueness of the  problem? 
With respect to question (ii), we prove in Theorem \ref{Theorem:PiecePhase} that if one is given an over-segmentation with \(M\) segments that satisfies \(N>2\tau+M,\) then this piecewise constant phase problem still admits a unique solution.   

Based on these theoretical results, in Section \ref{Section:SignAlgo} we address the first question above and develop methods to find either an exact or an approximate over-segmentation to intervals of constant sign, given only (noisy) measurements of $|\bf F|$. Given such an over-segmentation, we further develop a computationally efficient algorithm to retrieve the unknown sign pattern. Our approach follows our previous works \cite{OurIEEE,raz2011vectorial}, whereby instead of working with the unknown signal values ${\bf f}$ as our variables, we work with the unknown phases, and formulate for them a  quadratic functional to be minimized. Relaxing the requirement that the solution is a phase vector leads to solving a system of linear equations.
In the noise-free case we prove that with a sufficiently high sampling rate,  our algorithm is guaranteed to recover the true sign pattern.

In section \ref{Section:PR_Apps} we present two phase retrieval applications of practical interest where the sign problem arises.  The first is vectorial phase retrieval with only three measurement vectors. Here the problem is to recover two different compactly supported signals ${\bf f}_1$ and ${\bf f}_2$ from measurements of
$|{\bf F}_1|$, $|{\bf F}_2|$ and their interference $|{\bf F}_1+{\bf F}_2|$.
With sufficient over-sampling, this problem was recently proven to admit a unique solution in \cite{beinert2015ambiguities}, but no reconstruction algorithm
was given. The second application is 
the recovery of two well separated 1-D objects from a single spectrum, a problem known to have a unique solution \cite{crimmins1983uniqueness}. In this paper we show how stable recovery for both problems is possible by solving a related sign problem. Finally, in Section \ref{Section:Simulations} we illustrate the performance of our algorithm via several simulations. For an example with real 2-D experimental data (thus involving a 2-D sign problem), we refer the reader to \cite{leshem2016direct}. 

This discrete sign problem considered in this paper can be viewed as the finite dimensional analogue of a similar continuous problem recently studied by Thakur \cite{Thakur}. The problem considered in \cite{Thakur} is the recovery of a real-valued function $g(t)$ whose continuous Fourier transform $G(\omega)$ is band-limited, from discrete measurements \(|g(t_j)|\).
The analogy between the two problems follows by relating the continuous function $g$ to our finite dimensional vector $\bf F$ and its Fourier transform $G$ to $\bf f$. In \cite{Thakur}, Thakur proved that if the sampling rate is at least twice the Nyquist rate then the continuous sign problem is well posed. He further developed an algorithm to reconstruct the function from a finite number of \(N\) measurements. While with noise-free data the reconstruction error of Thakur's method decays exponentially fast in \(N\), as we illustrate in Section \ref{Section:Simulations}, it is quite sensitive to even small measurement noise. In contrast, our sign recovery algorithm, despite its underlying discrete formulation, is still applicable to this continuous setting, and offers improved robustness to noise. We conclude in Section \ref{Sec:discussion} with a discussion and directions for future research.

 %%%%%%%%%%%%%%%%%%%%%%%%%%%%%%%%%%%%%%%%
 %%%%%%%%%%%%%%%%%%%%%%%%%%%%%%%%%%%%%%%%
\section{Problem Setup}\label{Section:Setup}

%\textcolor{red}{First discuss \(f\in\mathbb{C}^n\) with compact support whose FT $F$ is real valued. 
%Given noisy measurements of $|F|$ we wish to recover the sign of $F$. 
%This problem raises several questions: (i) mathematical well-posedness: in the ideal case of noise-free discrete measurements, is it possible to uniquely recover the sign, up to a global \(\pm 1\) trivial sign ambiguity ?\  and (ii) is there a simple reconstruction algorithm to recover the sign, which is computationally efficient and robust to noise ?\ 
%Actually we solve a slightly different problem, where we assume $F$ has piecewise
%constant phase.} \\

 \textit{Notation.} We denote $\i=\sqrt{-1}$. For $z \in \mathbb{C}$ we denote by $\mathcal{R}(z)$ and $\mathcal{I}(z)$ its real and imaginary parts and by $\bar z$ its complex conjugate. Vectors appear in boldface letters, for example ${\bf F}=(F_0,\ldots,F_{N-1})$. We denote denote the entry-wise multiplication of two vectors $\bf F$ and $\bf G$ by ${\bf F}{\bf G}$. We further denote by ${\bf f}\star{\bf g}$ the cross-correlation between the vectors $\bf f$ and $\bf g$.

\textit{Measurements.} While in this paper we shall mostly study a discrete formulation, it is nonetheless instructive to first briefly review the continuous setting. Let $f_c(t)$ denote a one dimensional continuous signal, whose 1-D Fourier transform  $F_c(\omega)$ is given by
\begin{equation}
F_c(\omega)=\int f_c(t)e^{-\i\omega t}dt=|F_c(\omega)| e^{\i\phi(\omega)}.
\end{equation}
As mentioned in the introduction, in many applications direct measurement of $f_c(t)$ is not possible.  Rather, the  measured data is typically an equispaced sampling of $|F_c(\omega)|^2$ at $\omega_j=j\Delta\omega, \ j=0,...,N-1$. We denote the values of $F_c(\omega)$ at the sampled frequencies by $F_{j}:=F_c(\omega_j)$, and denote \({\bf F}=(F_0,\ldots,F_{N-1})\).  A discrete approximation of $f_c(t)$ at the \(N\)  points $t_k=\frac{2\pi k}{N\Delta\omega}$, \(k=0,\ldots,N-1\)
can be computed via the discrete Fourier transform ($\mathcal{DFT}$) of $F_{c}(\omega_j)$
\begin{equation}
f(t_k)=\frac{\Delta\omega}{2\pi}\sum_{j=0}^{N-1}|F_{c}(w_j)|e^{\i \phi(\omega_j)}e^{\i\omega_jt_k}:=
\mathcal{DFT}\left\{ {\bf F}\right\}
\label{Eq:DiscreteDef}
\end{equation} 
For convenience, we rescale the time and frequency coordinates so that $\Delta\omega = \frac{2\pi}{N}$ and $t_k=k$ for $k=0,...,N-1$. In this paper we restrict our attention to a finite dimensional formulation and consider the reconstruction of the \(N\) signal values ${\bf f}=(f_0,\ldots,f_{N-1})$ as our end goal. We note that, strictly speaking, $f(t_k) \neq f_c(t_k)$. However, for $N \gg 1$ and a sufficiently high sampling rate, $f(t_k) \approx f_c(t_k)$, see for example \cite{kammler2007first}.

For future use, we recall that the two vectors  $\bf f$ and ${\bf F}$ are thus related as follows
\begin{equation}
{f(k)}=\mathcal{DFT}\{{\bf F}\}(k)=\frac{1}{N}\sum_{j=0}^{N-1}F_je^{\i\omega_jk}
\ \mbox{and }
F_j=\mathcal{IDFT}\{{\bf f}\}(j)=\sum_{k=0}^{N-1}f(k)e^{-\i\omega_jk}.
\end{equation}
%\begin{align}
%& {f(k)}=\big(\mathcal{DFT}\{{\bf F}\}\big)(k)=\frac{1}{N}\sum_{j=0}^{N-1}F_je^{\i\omega_jk} %\\
%& F_j=\big(\mathcal{IDFT}\{{\bf f}\}\big)(k)=\sum_{j=0}^{N-1}f(k)e^{-\i\omega_jk}. %\nonumber
%\end{align}

\indent \textit{The classical phase problem.} Let $|{\bf F}|^2\in\mathbb{R}^N$ be the measured spectrum at the $N$ equispaced frequencies as described above and denote the unobserved phase
 vector by ${\bm \phi} = arg({\bf F})$. The phase problem is to reconstruct from the intensity measurements $|{\bf F}|^2$,  the missing phases $\left\{\phi_j\right\}_{j=0}^{N-1}$, or equivalently the \(N\) values ${\bf f}=\mathcal{DFT}\left\{{\bf F}\right\}$. Clearly, without additional constraints this problem is ill-posed, as any vector of phase values is a valid solution.
A commonly imposed constraint is that ${\bf f}$ has a compact support. However, even with this assumption, the 1-D phase retrieval problem  is in general still ambiguous. If $N>2\tau$ and \({\bf f}\) has a compact support of length \(\tau+1\), it still has at most $2^{\tau-1}$ different solutions, see   \cite{beinert2015ambiguities,Sodin}. 

\indent \textit{The sign problem.} In this paper we consider phase retrieval when the underlying function $F_c(\omega)$ is \textit{real-valued}. Namely, the unobserved phases of ${\bf F}$ are all restricted to the values $\phi_j=0$ or $\pi$.  Given $|{\bf F}|^2$ the problem is then to reconstruct this sign pattern.
  Clearly, without additional constraints this problem is also ill-posed as any sign pattern is a valid solution. In this paper we assume that ${\bf f}$ has a compact support and focus on the following two questions: 1) is there a unique solution? and 2) assuming there is, can one develop a stable and computationally efficient algorithm to recover it. 

\textit{Trivial ambiguities.} If \({\bf f}\) is a solution to the classical phase problem, then so are its circular shift, its reflection, conjugation and multiplication by a unimodular factor \(e^{{\bf i}\theta}\). These transformations do not change the physical nature of the signal and are thus known as trivial ambiguities of the phase problem. In the the sign problem, since $\bf F$ is real-valued it implies that $f(k)=f^*(-k \text{ mod } N)$ for all $k=0,\ldots,N-1$. This  eliminates the reflection ambiguity and all circular time shifts ambiguities except for a shift by $N/2$, assuming \(N\)\ is even. Hence, the trivial ambiguities of the sign problem are only a circular shift by $N/2$ and multiplication by a global sign.

\textit{Compactly supported signals.} We say that a signal ${\bf f} \in \mathbb{C}^N$ has a compact support of length $\tau +1$ if $f(k)=0$ for all $k \not\in CS$ where
\begin{equation}\label{Eq:CS_Phase}
CS=\left\{k:k \in [-\alpha,-\alpha+\tau]\mod N \right\} 
\end{equation}
for some $\alpha\in\mathbb{Z}$. For the sign problem, given the discussion above on trivial ambiguities and assuming for simplicity that $\tau$ is even,  there are only two options for the set $CS$. Either
\begin{equation}\label{Eq:CS_Sign}
CS=\left\{k:k\in[-\frac{\tau}{2},\frac{\tau}{2}] \mod N \right\}
\end{equation}
or its circular time shift by $N/2$. Without loss of generality, throughout this work we will use  Eq.\eqref{Eq:CS_Sign} as the assumed compact support of the sign problem.

\section{Mathematical Uniqueness}\label{Section:MathUniq}
\label{sec:unique}
%\textcolor{red}{Theorem of Uniqueness - goes through roots of high degree polynomials}

\noindent In this section we study the uniqueness of the sign problem.  First, we show that in the ideal setting of noise free measurements, if the signal $\bf f$ has  a sufficiently small compact support compared to $N $ (or equivalently, if $|F_c(\omega)|^2$ is sampled at a sufficiently high rate), then the sign problem admits a unique solution. Next, we show how the sign problem for compactly supported signals can be viewed as a special case of a more general phase problem, where the (complex-valued) phase is assumed to be piecewise constant. Moreover, we prove that at the expense of a higher sampling rate, this piecewise constant phase problem also admits a unique solution. Finally, in Section \ref{Section:SignAlgo} we combine the above results and develop a simple, computationally efficient algorithm to solve the original sign problem. All proofs appear in the appendix.
% Next, we prove uniqueness
%for the case in which the Fourier phase of a compact supported signal is piecewise constant, and a division to constant phase segments is known. Finally, we connect the former results by showing that the sign problem is actually a special case of a piecewise constant phase.
 
Specifically,  we introduce the following assumption:\ 
\begin{enumerate}
\item[A1] \label{Assump:CS}\textit{Compact support.} The vector $\mathbf{f}=\mathcal{DFT}\left\{ {\bf F}\right\}\in\mathbb{C}^N$ has a compact support of length \(\tau+1\). For the classical phase problem the compact support set is given by Eq.\eqref{Eq:CS_Phase}, whereas for the sign problem  it is given by Eq.\eqref{Eq:CS_Sign}.
\end{enumerate}

 In this section we assume that $\tau$ is known. In Section \ref{Section:SignAlgo} we discuss how $\tau$ can be estimated from the measured data. The following theorem proves that, when sampled properly, this sign problem has a unique solution.

 \begin{theorem}\label{Theorem:SignUniq}
  Let $|{\bf F}|^2$ be the observed noise-free spectrum of a signal that satisfies assumption A1. If ${\bf F}\in\mathbb{R}^N$ and $N>2\tau$ then the sign problem has a unique solution up to a global \(\pm 1\) sign ambiguity.
 \end{theorem}

 According to Theorem \ref{Theorem:SignUniq}, up to a $\pm 1$ global sign, there is a single sign pattern ${\bf s}=(s_0,\ldots,s_{N-1})$, such that $\mathcal{DFT}\{|{\bf F}|{\bf s}\}$ yields a signal with the correct compact support. As the following lemma shows, this sign pattern cannot be arbitrary. Rather, it has a limited number of sign changes.

 \begin{lemma}
                \label{lem:sign_changes}
  Let \({\bf F}\in\mathbb{R}^N\). If ${\bf f}=\mathcal{DFT}\{{\bf F}\}$ satisfies assumption A1, then the vector ${\bf s}=sign({\bf F}) \in \{-1,1\}^{N}$  has at most  $\tau$ sign changes. 
 \end{lemma}

% \begin{Remark} \label{Remark:F0}
% If $F(j)=0$ then ${\bf s}=sign({\bf F})$ is not well defined. \ Lemma \ref{lem:sign_changes} continues to hold, 
% \end{Remark}

 \begin{Remark} \label{Remark: SuperOss}
 It is interesting to contrast Lemma \ref{lem:sign_changes} with the phenomenon of {superoscillations}. As studied for example in \cite{ferreira2006superoscillations}, the Fourier transform of a compactly supported signal may have an arbitrarily large number of sign changes. As a specific example, the vector $F_j=(-1)^j$
has $N-1$ sign changes and its $\mathcal{DFT}$ is a delta function at the point $N/2$ which is obviously compactly supported. While at first sight this seems to contradict Lemma \ref{lem:sign_changes}, it does not, as Lemma \ref{lem:sign_changes} requires that the compact support is centered around the origin. Indeed, shifting the above delta function by $N/2$, so that its compact support is at the index \(k=0\) yields a  $\mathcal{DFT}$ vector \({\bf F}=(1,1,\ldots,1)\) with no sign changes at all, in accordance with Lemma \ref{lem:sign_changes}. 
 \end{Remark}
 
% \begin{Remark}\label{Remark:FiniteWindow}
% Lemma \ref{lem:sign_changes} is a consequence of our finite dimensional problem, which in particular implies that \({\bf F}\) has a finite number of sign changes. In a continuous setting,  $F_c(\omega)$ may have an infinite number of zeros. For example, if  $f_c(t)$ is the rectangular function, then $F_c(\omega)$ is the well known $sinc$ function, which has an infinite number of zeros. However, within a finite sampling window the number of zeros is finite.  Increasing the size of the window in which $F_c(\omega)$ is sampled may increase the number of zeros. This is equivalent to a finer sampling of $f_c(t)$ which results in a larger number of non-zero values of the finite vector $\bf f$. 
% \end{Remark}

  According to Lemma \ref{lem:sign_changes}, the sign pattern  \({\bf s}\)
that corresponds to a signal \({\bf f}\) with compact support of length $\tau+1$, belongs to a set of size $O(N^\tau)$. Unfortunately, the size of this set is exponential in \(\tau\). Thus, finding the unique solution to a given sign problem by exhaustive search over all possible sign patterns in this set is computationally intractable. 

To construct a computationally efficient sign retrieval algorithm, we  consider the following more general phase problem:\ We relax the strict assumption of a real-valued phase and instead assume the phase is complex-valued but piecewise constant in \(M\) a-priori known segments. Note that according to Lemma \ref{lem:sign_changes}, the (real-valued)\ sign problem is a particular instance of this piecewise constant phase problem. The following theorem shows that under suitable conditions, this modified phase problem also admits a unique solution.  In section \ref{Section:SignAlgo} below we then show how this modified problem can be solved computationally efficiently by framing it as the solution to a set of linear equations. 

We mathematically formulate the piecewise constant phase property as follows:\  
\begin{itemize}
\item[A2] \textit{Known segmentation into constant phase intervals:} There is a known division of the $N$ frequencies $\omega_j$ to $M$ contiguous segments, in which the unobserved phase vector ${\bm \phi}$ is piecewise constant. Let $\bf c$ be a vector of length $M$ containing the first index in each of these segments. The \(m\)-th segment of length $N_m$ consists of all frequencies $\omega_j$ with $j \in [c(m),c(m)+1,\ldots,c(m)+N_m-1]$ and $\phi(\omega_{c(m)})=\phi(\omega_{c(m)+k}), 1\leq k \leq N_m-1$.
%and is composed of the following frequencies, \(\omega_m,\omega_m+\Delta\omega,\ldots,\omega_m+(N_m-1)\Delta \omega\). Hence, the phase vector satisfies   for all \(m=1,\ldots,M \)
%\begin{equation*}
%e^{\i \phi(\omega_m+n\Delta \omega)}=e^{\i \phi(\omega_m+(n+1)\Delta \omega)}\quad n=0,...,N_m-2,
%\end{equation*}
%Clearly $\sum_{m=1}^M N_m=N$. 
\label{Assump:Seg} 
  \end{itemize}

\begin{Remark}\label{Remark:Seg}
Assumption A2 may define an over-segmentation of $\bm \phi$ to piecewise constant phase intervals. In particular, it does not imply that the number of intervals $M$ of the given segmentation  is minimal.
% That is, given ${\bf \phi}=arg({\bf F})$ which is composed of $\tilde M$ constant phase intervals, assumption \ref{Assump:Seg} may define a segmentation with $M>\tilde M$, but without assigning any indices in which $\bf \phi$ has different values to the same segment . In this case, we say that assumption \ref{Assump:Seg} defines an over-segmentation of $\phi_j$.
\end{Remark} 
 
 \begin{Remark}
If all indices where $\bf F$ changes its sign were exactly known (e.g., we had a complete rather than an over-segmentation), then this information would directly resolve the sign pattern \({\bf s}\), up to its inherent global \(\pm 1\) sign ambiguity. Unfortunately, in practice it is difficult to determine all zero crossings from the observed \(|{\bf F}|\), see for example \cite{Thakur}. Indeed, a key result of our paper, stated in Theorem  \ref{Theorem:PiecePhase} below, is that under suitable conditions even an over-segmentation can suffice to resolve the sign problem. \end{Remark}

%**** For notational convenience, we first define the shifted signal $F_{\alpha}(\omega_j)$  
% \begin{equation} \label{Eq:Shifted_F}
%F_{\alpha}(\omega_j)=e^{\i \omega_j \alpha} F(\omega_j)
% \end{equation}
%By the definition of the $\mathcal{DFT}$, the corresponding signal $f_\alpha(k)=f((k-\alpha)\mod N)$. In particular 
%for $\alpha=\tau/2$, from assumption (2), we obtain a shifted signal whose support is the set of indices $0,1,\ldots,\tau$.
%We denote this signal by $\mathbf{f}_s$ and its DFT by $F_s$. Note that $f_s(k)=\mathcal{DFT}\left\{F_s\right\}$ vanishes for $k\in[\tau+1,N-1]$.

 \begin{theorem} \label{Theorem:PiecePhase}
   Let $|{\bf F}|^2$ be the noise-free spectrum of a signal that satisfies assumptions A1 and A2. Assume $|F_j| \neq 0$ at the first and last indices of each of the \(M\) segments of the given segmentation. If $N>2\tau+M$ then the piecewise constant phase problem has a unique solution up to multiplication by a global phase.
\end{theorem}

\begin{Remark}
Our proof of Theorem \ref{Theorem:PiecePhase} breaks down if \(F_j=0\) at an end-point of one of the segments. However, the proof can be modified to cover this case, as well as the case of a number of consecutive points in which $|{\bf F}|$ vanishes near segment edges. See the appendix for details.  
\end{Remark}

\section{Sign reconstruction with known support}\label{Section:SignAlgo}
We now describe a simple and computationally efficient algorithm to solve the sign problem. First, we study the case of noise-free measurements and further assume that an over-segmentation of \([0,1,\ldots,N-1]\) into intervals of constant sign is a-priori known. Relying on Theorem 2, we show in Section
\ref{SubSec:Known_Seg} how the sign pattern can then be cast as the solution to an over-determined system of linear equations. 
Next, in Section \ref{SubSec:Seg} we describe a method to segment the set of \(N\) indices into intervals of constant sign, given only the vector \(|{\bf F}|\). An algorithm to recover the sign pattern in the presence of noisy measurements appears in Section \ref{SubSec:SR_Alg}. Finally, the estimation of the typically unknown compact support length \(\tau\) is addressed in Section \ref{SubSec:tau_estimation}. 
%Lemma 1 shows that when $f(k)$ has a compact support the sign of $F(\omega_j)$ is not arbitrary but rather segmented to constant sign regions. In this case, as we detail below, Theorem 2 suggests a simple algorithm, in which the phase constraint is relaxed. 

%%%%%%%%%%%%%%%%%%%%%%%%%%%%%%%
\subsection{Sign recovery with known over-segmentation}
        \label{SubSec:Known_Seg}

Consider a signal ${\bf F}\in\mathbb{R}^N$ satisfying assumptions A1 and A2, and further assume that the length \(\tau\) of the compact support of  $\bf f$ is a-priori known. As in our earlier works \cite{DBFH,OurIEEE}, instead of working with the \((\tau+1)\) signal values \(f(k)\), $k=-\tau/2,\ldots,\tau/2\, (\!\!\!\!\mod N)$, as the unknown variables, we consider the vector of $N$ unknown signs ${\bf X}=(X_0,\ldots,X_{N-1})$. Further, rather than dealing with the combinatorial set \({\bf X}\in\{-1,1\}^N\), we relax this constraint and allow all entries \(X_{j}\) to be complex valued, namely  ${\bf X}\in\mathbb{C}^N$. 

As we now show, this allows us to write a system of linear equations for the vector ${\bf X}$ over the field $\mathbb{C}$ whose unique real-valued solution is the true sign pattern. The first set of equations captures the compact support assumption on the vector \({\bf f}={\mathcal DFT}\{{\bf F}\}\). Since ${\bf f}$ is zero outside the set of indices CS of Eq. (\ref{Eq:CS_Sign}), the unknown vector ${\bf X}$ must satisfy the following set of \(N-\tau-1\) linear equations
\begin{equation}
        \label{Eq:X_CS}
\mathcal{DFT}\left\{|\bf F|{\bf X}\right\}(k)=0 \ , \ k \notin CS.
\end{equation} 
The second set of equations imposes the known segmentation into intervals of constant sign, as described in A2. Let $\mathcal{S}$ denote the set of indices within the intervals of constant sign, excluding the last one in each interval. Then,  
\begin{equation}
        \label{Eq:X_ConstPhase}
X_{j}={X}_{j+1} ,  \ \forall j\in\mathcal{S}. 
\end{equation}
Given a segmentation to $M$ constant sign intervals the number of equations in (\ref{Eq:X_ConstPhase}) is thus \(N-M\). 

By definition, the true sign pattern \({\bf s}=sign({\bf F})\) is a solution of Eqs. (\ref{Eq:X_CS})-(\ref{Eq:X_ConstPhase}). Theorem \ref{Theorem:SignAlgo} below shows that under suitable conditions and up to multiplication by a constant, it is the \textit{only}  solution to this system of equations. 

\begin{theorem} \label{Theorem:SignAlgo}
Let $|{\bf F}|^2$ be the intensity of a signal that satisfies assumptions A1 and A2.   Assume $|F_j| \neq 0$ at the first and last indices of each of the \(M\) segments of the given segmentation. If $N>2\tau+M$ then, in the noise-free case, the only solutions to the set of linear equations (\ref{Eq:X_CS})-(\ref{Eq:X_ConstPhase}) are of the form \({\bf X}=c\,{\bf s}\) where $c\in\mathbb{C}$. \end{theorem}

\begin{Remark}
Importantly, any over-segmentation of $\bf X$ with number of segments $M<N-2\tau$ suffices for Theorem \ref{Theorem:SignAlgo} to hold. That is, not all the indices for which the true sign pattern $\bf s$ is piecewise constant need to be captured in Eq. \eqref{Eq:X_ConstPhase}. 
\end{Remark}

\begin{Remark}
While the condition \(N>2\tau+M\) is sufficient to ensure a rank-one solution to the linear system (\ref{Eq:X_CS})-(\ref{Eq:X_ConstPhase}), it is by no means a necessary condition. Empirically, often an over-segmentation with more segments still suffices to reconstruct the correct sign pattern. 
\end{Remark}

\subsection{Segmentation to constant sign intervals} \label{SubSec:Seg}

By Theorem \ref{Theorem:SignAlgo}, with a suitable over-segmentation one can retrieve the sign of  $\bf F$ by solving the set of linear equations (\ref{Eq:X_CS})-(\ref{Eq:X_ConstPhase}). To this end, however, one must first determine an (over-)segmentation from $|{\bf F}|^2$ alone. The following lemma provides a principled method to do so. 

\begin{lemma}\label{lem:CrossTh}
        Let ${\bf F}\in \mathbb{R}^N$ be a signal that satisfies A2. Then, the difference between two consecutive values of $\bf F$ is bounded by
        \begin{equation}\label{Eq:LemmaBound}
        |F_j-F_{j-1}| \leq  \big(\frac{2}{N}\big)^{3/2}\pi {\mathcal S_{\tau/2}} \|{\bf F}\|
        \end{equation}
where $\|{\bf F}\|^2= \sum_{j=0}^{N-1} F_j^2$ and 
\begin{equation*}
\mathcal S_{\tau/2}^2 = \sum_{k=1}^{\tau/2} k^2=\frac{\tau(\tau+1)(\tau+2)}{24}=O(\tau^3). 
\end{equation*}
\end{lemma}

Lemma \ref{lem:CrossTh} implies that at any pair of consecutive indices $\{j-1,j\}$ that satisfy
\begin{equation}\label{Eq:BoundCond} 
|F_j|+|F_{j-1}|>\left(\frac{2}{N}\right)^{3/2}\pi {\mathcal S_{\tau/2}} \|{\bf F}\|
\end{equation}
the sign of $\bf F$ must be equal. Hence, knowing only the vector $|{\bf F}|$, Eq. \eqref{Eq:BoundCond} provides an over-segmentation to constant sign intervals. If the number of found segments $M$ is sufficiently small so that $N>2\tau+M$ then we can directly recover the sign pattern by solving the system of linear equations \eqref{Eq:X_CS}-\eqref{Eq:X_ConstPhase}.

\begin{Remark} \label{Remark:Large_N}
With a sufficiently high oversampling rate, the number of segments determined by Eq. (\ref{Eq:BoundCond}) satisfies \(N>2\tau+M\), which in turn guarantees recovery of the true sign pattern by  solving equations  \eqref{Eq:X_CS}-\eqref{Eq:X_ConstPhase}. To see this, note that as we increase the number of measurements \(N\) while keeping \(\tau\) fixed,  \(\|{\bf F}\|\) increases as \(O(\sqrt{N})\) and hence the threshold on the right hand side of Eq. (\ref{Eq:BoundCond}) decreases as \(O(1/N)\). Thus, for sufficiently large \(N\), a sufficient number of pairs \((j-1,j)\) satisfy Eq. (\ref{Eq:BoundCond}). 
\end{Remark}

Unfortunately, even though  Eq. \eqref{Eq:BoundCond} provides a correct over-segmentation, with a finite over-sampling rate the resulting number of segments can be too large, so that the corresponding Eqs. \eqref{Eq:X_CS}-\eqref{Eq:X_ConstPhase} have multiple solutions. To decrease the number of segments, and thus increase the number of linear equations, we supplement Eq. \eqref{Eq:BoundCond} with a heuristic segmentation scheme which works very well in practice, even though it is not theoretically guaranteed to yield a correct over-segmentation. Our heuristic segmentation scheme relies on the fact that for a sign change to occur between ${F}_j$ and ${F}_{j+1}$ its analytic continuous extension $F(\omega)$ must have a zero crossing at some intermediate $\omega \in [\omega_j,\omega_{j+1}]$. Hence, $|{\bf F}|$ is likely to have a local discrete minimum near this zero crossing. 

This observation leads to the following segmentation algorithm summarized in Table \ref{Tab:HeuSeg}:  First, the discrete local minima of $|{\bf F}|$ are found and their indices are defined as single index segments. Next, all the indices between each pair of consecutive minima are defined to have the same sign. Finally, for each minimum index, the adjacent index for which  $|{\bf F}|$ has closer value to the minimum value is also defined as a single index segment. The last step is applied to reduce errors in the resulting segmentation. In the absence of noise, our final segmentation is the merging of both segmentations described above. Fig. \ref{Fig:Seg_Example} in Section \ref{SubSec:SegScheme} presents an example of the segmentations produced by both schemes. 
%The indices of the true sign changes are denoted by blue stars. The indices defined by Eq. \eqref{Eq:BoundCond} are marked by red circles, illustrating the resulting over-segmentation. The indices marked by green circles denote the additional indices defined by our heuristic segmentation scheme according to Table \ref{Tab:HeuSeg}, decreasing the overall number of segments.  

\begin{table}[H]
\begin{center}
\begin{tabular}{l}
\hline
{\bf Algorithm} Heuristic over-segmentation\\ \hline\hline\\
{\bf Input: } $|{\bf F}|^2$.\\
{\bf Algorithm: }\\
1: Find the local minima of  $|{\bf F}|$ given by the indices $j$ s.t. $|{F}_j| <  \min\{|F_{j-1}|,|F_{j+1}|\}$.  \\
2. Define all minima indices as single index segments. \\
3. Define the indices between each pair of consecutive minima indices as \\
 a constant sign interval. \\
4: For each minima index, find the adjacent index with the closest $|{\bf F}|^2$ \\
value and exclude it from the constant sign interval, i.e. define it as a single index segment.\\
{\bf Output: } Over-segmentation to constant sign intervals. \\
\hline
\end{tabular}
\end{center}
\caption{Heuristic over-segmentation scheme.}
\label{Tab:HeuSeg}
\end{table}

\subsection{Sign retrieval from noisy measurements}
        \label{SubSec:SR_Alg} 
    In any realistic scenario the vector of intensities  $|{\bf F}|^2$ is measured with some noise. In this case, no sign pattern satisfies Eq. \eqref{Eq:X_CS}. To cope with measurement noise, we reformulate the sign problem as the minimization of a suitable quadratic functional. To this end, we first describe our noise model, then construct the functional to be minimized and detail our minimization approach.
   
   Let $|\tilde {\bf F}|^2$ be the vector of noisy measurements. As described in \cite{OurIEEE}, a rather general noise model is    \begin{equation}\label{Eq:NoiseModel}
  |\tilde {\bf F}|^2 = |{\bf F}+{\tfrac{\sigma}{\sqrt{N}}\bm \eta}^s|^2+|{\bf F}+{\tfrac{\sigma}{\sqrt{N}}}{\bm \eta}^s|{\bm \eta}^{sh}+{\bm \eta}^d
   \end{equation}
   where ${\bm \eta}^s$ is background additive noise with noise level \(\sigma\), ${\bm \eta}^d$ is dark counts noise, and $|{\bf F}+{\tfrac{\sigma}{\sqrt{N}}}{\bm \eta}^s|{\bm \eta}^{sh}$ is the detector shot noise which is proportional to the signal intensity. In many cases the dark counts ${\bm \eta}^d$ has small variance and its main effect can be removed by proper calibration. We thus assume for simplicity that ${\bm \eta}^d=0$.
We assume a classical light experiment, far from the single photon regime and we thus neglect the effect of shot noise. We further assume that ${\bm \eta}^s=(\eta^s_0, \ldots , \eta^s_{N-1})$ consists of $N$ independent and identically distributed (i.i.d.) complex-valued Gaussian random variables $\mathcal{N}(0,1)$.  

To cope with noisy measurements, we modify our sign retrieval scheme in two ways:\ (i) refine the method to find an oversegmentation; (ii) replace the solution of a homogeneous linear system by the minimization of a quadratic functional. 
In the presence of noise, our sign retrieval consists of the following key steps:

\noindent {\bf Step 1: Find an over-segmentation}.
 First we apply the heursitic segmentation scheme defined in Table \ref{Tab:HeuSeg}. We denote by  $\mathcal{S}_1$ the set of indices of size $S_1$ in the resulting  segmentation, excluding the last index in each segment. Next, we modify Eq. \eqref{Eq:BoundCond} as detailed in the Appendix to yield,
\begin{equation}\label{Eq:BoundCond_Noisy} 
|F_j|+|F_{j-1}|>\left(\frac{2}{N}\right)^{3/2}\pi {\mathcal S_{\tau/2}} \|\tilde{\bf F}\|+\sigma/\sqrt{N}
\end{equation}
where $\sigma$ is the standard deviation of the noise.  We denote by  $\mathcal{S}_2$ the set of indices of size $S_2$ in the segmentation defined by Eq. \eqref{Eq:BoundCond_Noisy}, excluding the last index in each segment.
   
   \noindent {\bf Step 2: Construct and minimize a quadratic functional}.
   Given the noisy measurements $|\tilde {\bf F}|^2$ and the compact support parameter $\tau$, we define the matrix $\mathcal{A}_{cs}$ of size $(N-\tau-1 )\times N$ as
   \begin{equation}\label{Eq:A_Def}
      (\mathcal{A}_{cs})_{k,j} = (DFT \cdot \mathcal{D})_{k,j}, \  \ k \notin CS, \ j\in[1,N]
   \end{equation}
   where $(DFT)_{j,k}=\frac{1}{N}e^{\i\omega_{j-1}(k-1)}$ with $\ j,k=0,\ldots,N-1$ is the discrete Fourier transform matrix, $\mathcal{D}=diag(|\tilde {\bf F}|)$ is a diagonal matrix and   the set $CS$ is given by Eq. \eqref{Eq:CS_Sign}. In the noise free case, $\mathcal{A}_{cs}X=0$ imposes that $\mathcal{DFT}\{|\tilde {\bf F}|{\bf X}\}$ vanishes for $k \notin CS$ and is precisely Eq. \eqref{Eq:X_CS}. 
   %In the case of noisy measurements instead seek to minimize $\|\mathcal{A}_{cs}\bf{X}\|^2$.
   
   Next, we impose our two segmentation schemes starting with the heuristic one. We use $\mathcal{S}_1$ to  construct the matrix $\mathcal{A}_1$ of size $S_1 \times N$, whose $k^{th}$ row is given by
   \begin{equation}\label{Eq:S1_Mtx}
   (\mathcal{A}_1)_{k,j}= 
   \begin{cases}
    W(k) & j=l_k \\
   -W(k) & j=l_k+1 \\
   0 & \text{otherwise.}
   \end{cases}  \ \ \ \ l_k \in \mathcal{S}_1 
   \end{equation}
    where $k\in[1,S_1]$, $l_k$ is the $k$-th index in \({\cal S}_1\), and the vector of weights ${\bf W}=[W_1,\dots,W_{S_1}]$ is defined as $W(k)=\min\{|{\tilde F}_{l_k}|^2,|{\tilde F}_{l_k+1}|^2\}$. We note that $\mathcal{A}_1{\bf X}=0$ imposes that $\bf X$ is piecewise constant in the intervals defined by Table \ref{Tab:HeuSeg}. The purpose of $\bf W$ is to account for the fact that the heuristic algorithm presented in Table \ref{Tab:HeuSeg} does not guarantee a correct over-segmentation and that errors typically occur at indices corresponding to low $|\tilde {\bf F}|^2$ values where we have less certainty in our segmentation. With our weighting approach, we have that $\mathcal{A}_1{\bf X}=0$ imposes  $W_j(X_j-X_{j+1})=0$ for $j \in \mathcal{S}_1$. Hence, when a quadratic functional based on $\mathcal{A}_1$ is minimized, as described below, the constraint that $\bf X$ is piecewise constant at $j\in \mathcal{S}_1$ is suppressed for entries corresponding to low $|\tilde {\bf F}|^2$ values. 
    
    Our second segmentation approach is guaranteed to yield a correct over-segmentation in the noise-free case, and in practice for high SNR, it is unlikely to make errors. Hence, we impose this segmentation without weights, in a way that also reduces the computational cost of our problem by reducing the number of variables. To this end, we sum over the appropriate columns of the $(N+S_1-\tau-1) \times N$ matrix $[\mathcal{A}_{cs};\mathcal{A}_1]$ to construct the matrix $\mathcal{A}$ according to    
    \begin{equation}\label{Eq:Def_A}
         \mathcal{A}_{k,m}=\sum_{j\in \text{segment m}}([\mathcal{A}_{cs};\mathcal{A}_1])_{k,j}, \ 1\le m \le M, \ k\in [1, N+S_1-\tau-1]. 
    \end{equation}  
%    \begin{equation}\label{Eq:Def_A}
%     \mathcal{A}_{k,j}=([\mathcal{A}_{cs};\mathcal{A}_1])_{k,j}+([\mathcal{A}_{cs};\mathcal{A}_1])_{k,j+1}, \ j\in \mathcal{S}_2, \ k\in [1, N+S_1-\tau-1]. 
%     \end{equation} 
    The size of $\mathcal{A}$ is $(N+S_1-\tau-1) \times (N-S_2)$. 
     
  Our approach is to find a vector $\bf X$ which minimizes the following,
\begin{equation}\label{Eq:Qx1}
        Q({\bf X}) = \|\mathcal{A}{\bf X}\|^2.
\end{equation}
Since the required output is a sign vector, in principle Eq. \eqref{Eq:Qx1} should be minimized over the $\{-1,1\}^N$, which results in a non-convex problem. Here we relax this constraint and instead minimize $Q({\bf X})$ over the set of complex-valued vectors. Similarly to \cite{OurIEEE}, given the $\pm 1$ global sign ambiguity, without loss of generality, we set $X(\nu)=-1$ for some index $\nu \in \{0, \ldots ,N-1\}$. The functional $Q({\bf X})$ of Eq. \eqref{Eq:Qx1} can thus be written as
   \begin{equation}\label{Eq:Qx2}
   Q({\bf X}) =  \|\mathcal{A}_{-\nu}{\bf X}_{-\nu}-a_{\nu}\|^2
   \end{equation}
   where $a_{\nu}$ is the $\nu^{\text{th}}$ column of $\mathcal{A}$, $\mathcal{A}_{-\nu}$ is the matrix $\mathcal{A}$ without $a_{\nu}$ and ${\bf X}_{-\nu}$ is the vector $\bf X$ without its $\nu^{\text{th}}$ entry. We choose $\nu$ as the index in which $|\tilde {\bf F}|^2$ is maximal.
   Minimizing $Q({\bf X})$ amounts to the convex problem of solving a set of linear equations without constraints. Moreover, in the noise-free case, minimizing $Q({\bf X})$ is equivalent to solving the set of linear equations
      \begin{equation}\label{Eq:SignAlgoMtx}
      \mathcal{A}{\bf X}=0
      \end{equation}
      which effectively imposes Eqs. \eqref{Eq:X_CS}-\eqref{Eq:X_ConstPhase} with $\mathcal{S} = \mathcal{S}_1 \cup \mathcal{S}_2$. 
       Therefore, under the assumptions of Theorem \ref{Theorem:SignAlgo} the minimizer of $Q({\bf X})$ is the true sign vector $\bf s$ (up to multiplication by a constant). As demonstrated in section \ref{SubSec:Noisy_Sim}, empirically this approach is robust to noise and, in practice, even tolerates a small number of segmentation errors. Table \ref{Tab:SignRet} summarizes our algorithm.\\
\begin{table}[H]
\begin{center}
\begin{tabular}{l}
\hline
{\bf Algorithm} Sign retrieval\\ \hline\hline\\
{\bf Input: } $|\tilde {\bf F}|^2$, $\tau$, $\sigma^2$.\\
{\bf Algorithm: }\\
1: Compute the sets of indices $\mathcal{S}_1$ and $\mathcal{S}_2$ according \\
\ \ \ \ to Table \ref{Tab:HeuSeg} and Eq. \eqref{Eq:BoundCond_Noisy}  respectively.\\
2: Construct the matrix $\mathcal{A}$. \\
3. Compute the minimizer $\hat {\bf X}$ of $\|\mathcal{A}_{-\nu}-{\bf X}_{-\nu}-a_{\nu}\|^2$.\\
4. Project $\hat {\bf X}$ onto a sign: $\hat {\bf s} = sign({\mathcal{R}(\hat {\bf X}}))$. \\
{\bf Output: } Estimated sign vector: $\hat {\bf s}$ \\
\hline
\end{tabular}
\end{center}
\caption{Sign retrieval algorithm.}
\label{Tab:SignRet}
\end{table}   

\subsection{Estimating the compact support}
\label{SubSec:tau_estimation}

In practice, the length $\tau+1$ of the compact support of $\bf f$ is typically not known precisely. 
Here we propose an algorithm to estimate it from the (possibly noisy)\ measurements $|\tilde {\bf F}|^2$, under the assumption that $\tau_{min} \leq \tau \leq\tau_{max}$ for some a-priori known values of $\tau_{min}$ and $\tau_{max}$. 

To estimate \(\tau\), we scan over the possible values   $\tau_{min} \leq \tau_s \leq\tau_{max}$. For each value  $\tau_s$ we retrieve the sign pattern $\hat {\bf s}$ according to the algorithm in Table \ref{Tab:SignRet}, compute the corresponding signal $\hat {\bf f}=\mathcal{DFT}\{|{\bf F}|\hat{\bf s}\}$ and its average energy outside of the assumed compact support, 
\begin{equation}
E_{out}(\tau_s)=\frac{1}{N-\tau_s-1}\sum_{k\notin CS(\tau_s)} |\hat f(k)|^2.
\end{equation} 
Our estimate $\hat \tau$ of \(\tau\) is the first location where \(E_{out}(\tau_s)\) attains its minimal value. This simple scheme is summarized in Table \ref{Tab:EstCS}.

To justify this approach let us first analyze the noise-free case. Here, for $\tau_s=\tau$, as proven in Section \ref{SubSec:Known_Seg} our algorithm perfectly recovers the true signal $\bf f$. Hence, at the correct compact support, $E_{out}(\tau)=0$. 
For $\tau_s<\tau$ there is no signal \({\bf f}\) with compact support parameter \(\tau_s\) whose DFT\ has a piecewise constant-phase. Namely, there is no vector \({\bf X}\) which gives \(Q({\bf X})=0\) in Eq. (\ref{Eq:Qx1}). Minimizing this quadratic functional gives as output some signal ${\bf f}$ which does not vanish outside the assumed compact support, as otherwise this would contradict Theorem \ref{Theorem:SignAlgo}. 
Hence, for any $\tau_s<\tau$, $E_{out}(\tau_s)>0$.
For $\tau_s>\tau$, the only solution is the true sign vector ${\bf s}$. In the noise-free case \(\hat \tau = \tau\).
  
In the presence of noise, even at \(\tau_s=\tau\), due to the noise, the recovered signal is noisy as well, and \(E_{out}(\tau)>0\). As $\tau_s$ increases above $\tau$ the number of linear equations in Eq. \eqref{Eq:SignAlgoMtx} decreases. Hence, the sensitivity to noise of the corresponding solution is increased, which in turn leads to an increase in  $E_{out}(\tau_s)$  as a function of $\tau_s$. When the noise level is low, our approach is thus still able to correctly estimate the true support parameter.

\begin{table}[H]
\begin{center}
\begin{tabular}{l}
\hline
{\bf Algorithm} Compact support estimation\\ \hline\hline\\
{\bf Input: } $|\tilde {\bf F}|^2$, $\tau_{min}$, $\tau_{max}$, $\sigma^2$.\\
{\bf Algorithm: }\\
1: For $\tau_s = \tau_{min}$ to $\tau_{max}$ do:  \\
\ \ \ \ \ \ \ a: Retrieve the sign of $\bf F$ using Table \ref{Tab:SignRet} with compact support of $\tau_s +1$. \\
\ \ \ \ \ \ \ b: For the retrieved sign $\hat {\bf s}$, compute $\hat {\bf f}=\mathcal{DFT}\{|\tilde {\bf F}|\hat {\bf s}\}$.\\
\ \ \ \ \ \ \ c: Compute $E_{out}(\tau_s)=\sum_{k\notin CS} |\hat f(k)|^2 / (N-\tau_s-1)$. \\
2. Estimate $\tau$ by $\hat \tau=\underset{\tau_s}{\text{argmin}} \{E_{out}(\tau_s)\}$. \\
{\bf Output: } Estimated compact support: $\hat \tau+1$ \\
\hline
\end{tabular}
\end{center}
\caption{Compact support estimation.}
\label{Tab:EstCS}
\end{table}

\section{Phase retrieval applications of the sign problem}\label{Section:PR_Apps}

We now present two phase retrieval settings of practical interest in which the sign problem plays a key role. The first is vectorial phase retrieval (VPR) \cite{DBFH,OurIEEE} but with only 3 measurements, and the second is phase retrieval from two sufficiently separated objects \cite{crimmins1983uniqueness,leshem2016direct}. 

 \subsection{VPR with 3 measurements} \label{SubSec:VPR_3M}
 VPR consists of a recently suggested family of physically feasible measurement schemes together with computationally efficient methods to recover the phase. In VPR, the problem is to recover two compactly supported signals ${\bf f}_1, {\bf f}_2 \in \mathbb{C}^N$ with corresponding Fourier transforms ${\bf F}_1$ and ${\bf F}_2$  from the following measurements,  
\begin{align}\label{Eq:VPRinput}
|{\bf F}_1|, \quad |{\bf F}_2|, \quad {\bf F}_1{\bf F}^*_2.
%& {\bf E}_1 = |{\bf F}_1| \\ \notag
%& {\bf E}_2 = |{\bf F}_2| \\ \notag
%& {\bf E}_3 = {\bf F}_1{\bf F}^*_2
\end{align} 
 As proven in \cite{OurIEEE},  under suitable conditions this phase problem admits a unique solution. Furthermore, the phase vectors  ${\bf X}_1, {\bf X}_2 \in \mathbb{C}^N$ that correspond to ${\bf F}_1$ and ${\bf F}_2$  can be uniquely retrieved by solving the following set of linear equations   
 \begin{equation}
 \begin{aligned} \label{Eq:VPReqs}
 & \mathcal{DFT}\{|{\bf F}_1|{\bf X}_1\}(k)=0, \ k\notin CS \\ 
 &\mathcal{DFT}\{|{\bf F}_2|{\bf X}_2\}(k)=0, \ k\notin CS  \\ 
 & {\bf F}_1{\bf F}^*_2{\bf X}_2=|{\bf F}_1||{\bf F}_2||{\bf X}_1.
 \end{aligned}
 \end{equation}

 In \cite{OurIEEE}, several physical scenarios were described where the following 4 vectors can be measured, $|{\bf F}_1|$, $|{\bf F}_2|$, $|{\bf F}_1+{\bf F}_2|$ and $|{\bf F}_1+\i {\bf F}_2|$. From these measurements the interference term ${\bf F}_1{\bf F}^*_2$ can be easily computed as $\frac{1}{2}(\big|{\bf F}_1+{\bf F}_2|^2+\i|{\bf F}_1+\i {\bf F}_2|^2-(1+\i)(|{\bf F}_1|^2+|{\bf F}_2|^2)\big)$. Then, the phase is retrieved by solving  Eq.\eqref{Eq:VPReqs}.  However, in various physical scenarios obtaining the fourth measurement, $|{\bf F}_1+\i {\bf F}_2|$, is difficult or impossible and only the following three spectra can be measured 
    \begin{align}\label{Eq:VPR3}
   |{\bf F}_1|, \quad |{\bf F}_2|, \quad {\bf S} = |{\bf F}_1+{\bf F}_2|.
   \end{align}
   As we now show, using sign retrieval, these three measurements allow to recover  ${\bf F}_1{\bf F}^*_2$, the interference term required to apply VPR to solve the phase problem.
   
  It was recently proven in \cite{beinert2015ambiguities}, that the phase problem corresponding to Eq. \eqref{Eq:VPR3} admits a unique solution if ${\bf f}_1$ and ${\bf f}_2$ have a sufficiently small compact support, but without suggesting a possible algorithm to solve it. 
We propose the following scheme: Given  Eq.\eqref{Eq:VPR3}, compute
   \begin{equation} \label{Eq:Cos}
   {\bf E}_R=\frac{1}{2}({\bf S}^2-|{\bf F}_1|^2-|{\bf F}_2|^2)=|{\bf F}_1||{\bf F}_2|\cos(\phi_{12})= \mathcal{R}({\bf F}_1{\bf F}^*_2)
   \end{equation}
where ${\bm \phi}_{12}=arg({\bf F}_1{\bf F}_2^*)$. Since ${\bf E}_R=\mathcal{R}\big[{\bf F}_1{\bf F}_2^* \big]$ its  $\mathcal{DFT}$ is compactly supported. This is also true for the unknown corresponding imaginary part  ${\bf E}_I=\mathcal{I}\big[{\bf F}_1{\bf F}_2^* \big]$. Moreover, the absolute value of ${\bf E}_I$ can be calculated from Eq.\eqref{Eq:VPR3},
\begin{equation*}
|{\bf E}_I|=\sqrt{|{\bf F}_1|^2|{\bf F}_2|^2-{\bf E}_R^2}.
\end{equation*}
To fully recover ${\bf F}_1{\bf F}_2^*$ we thus need to solve the following sign problem: Given the absolute value of the imaginary part of the interference term  ${\bf E}_I$, whose $\mathcal{DFT}$ is compactly supported, retrieve $sign({\bf E}_I)$. We can thus use our sign retrieval algorithm to retrieve this sign.  This immediately allows computation of ${\bf E}_R+\i{\bf E}_I={\bf F}_1{\bf F}_2^*$ which, together with Eq. \eqref{Eq:VPR3} completes the required input for VPR depicted in Eq.\eqref{Eq:VPRinput}. Applying this scheme, we demonstrate VPR reconstructions with 3 measurements in Section \ref{SubSec:VPR_3M_Sim}.

\subsection{Phase retrieval from separated objects}\label{SubSec:VPR_SepObjs}
A second phase retrieval scenario in which the sign problem arises is the reconstruction of two compactly supported objects that are well-separated, by more than the length of the larger compact support. Let ${\bf f}\in \mathbb{C}^N$ be of the form ${\bf f}={\bf f}_1+{\bf f}_2$ where ${\bf f}_1$ and ${\bf f}_2$ are compactly supported and well-separated. In \cite{crimmins1983uniqueness}, it was shown that despite this being a 1-D phase problem, under suitable conditions, the signal ${\bf f}$ is uniquely determined by the single measurement vector $|{\bf F}|=|\mathcal{IDFT}\{{\bf f}\}|=|{\bf F}_1+{\bf F}_{2}|$.  

Given the above uniqueness result, the goal here is to reconstruct $\bf f$ from the single measurement $|{\bf F}|^2$. Our approach is to use sign retrieval as a step to recover the input required for VPR from $|{\bf F}|^2$ and then apply VPR to retrieve $\bf f$. Since ${\bf f}_1$ and ${\bf f}_2$ are well-separated, we have that $\mathcal{DFT}\left\{|{\bf F}_1+{\bf F}_2|^2\right\}$ yields 3 separated terms. The central term is the sum of the autocorrelations of ${\bf f}_1$ and ${\bf f}_2$, given by ${\bf f}_1 \star {\bf f}_1+{\bf f}_2 \star {\bf f}_2$. The two other terms are their cross-correlation, ${\bf f}_1 \star {\bf f}_2$ and its complex conjugate. Performing an inverse $\mathcal{DFT}$ separately on each of these terms, yields the following equations
\begin{align} \label{Eq:SepObj}
{\bf I}_S = |{\bf F}_1|^2+|{\bf F}_2|^2, \quad {\bf E}_3 = {\bf F}_1{\bf F}_2^*. 
\end{align}
In order to apply VPR we first need to resolve the vectors $|{\bf F}_1|^2$ and $|{\bf F}_2|^2$ from Eq. \eqref{Eq:SepObj}. To this end, consider the (unknown) function 
\begin{equation}
{\bf I}_D=|{\bf F}_1|^2-|{\bf F}_2|^2
\end{equation}
Since ${\bf I}_D$ is the difference between the Fourier intensities of ${\bf f}_1$ and ${\bf f}_2$, its Fourier transform is also compactly supported. Furthermore, its absolute value can be computed from Eq.\eqref{Eq:SepObj} as,
\begin{equation*}
|{\bf I}_D|=\sqrt{{\bf I}_S^2-4|{\bf E}_3|^2}.
\end{equation*}
This gives rise to the following sign problem: Given $|{\bf I}_D|$ with $\mathcal{DFT}\{{\bf I}_D\}$  compactly supported, retrieve $sign({\bf I}_D)$. Once the sign is retrieved $|{\bf F}_1|^2$ and $|{\bf F}_2|^2$ can be immediately computed from their sum and difference which together with ${\bf F}_1{\bf F}_2^*$ gives the required input of Eq. \eqref{Eq:VPRinput}. The underlying two signals $\bf{f}_1$ and ${\bf f}_2$ can be now recovered by applying VPR.

 In Section \ref{SubSec:VPR_SepObjs_Sim} below we numerically demonstrate this scheme. For an application of this scheme to the  reconstruction of separated 2D objects from experimental X-ray free electron laser measurements, see \cite{leshem2016direct}.  

\section{Simulations}\label{Section:Simulations}

We illustrate the performance of our algorithm via several simulations. 
First, we consider noise-free and noisy reconstructions of complex-valued, compactly supported random signals. Next, we consider the continuous sign problem, and compare our method to the one suggested by Thakur \cite{Thakur}. Finally, we apply our sign retrieval algorithm  to the two phase retrieval applications described in Section \ref{Section:PR_Apps}, VPR with 3 measurements  and phase retrieval for two well-separated 1-D objects.

Given the global \(\pm1\) ambiguity of the sign problem, we measure the reconstruction quality by the  following mean-square-error (MSE),
\begin{equation}
\text{MSE}=\min\big(\sum_k |{\hat f(k)}-f(k)|^2, \ \sum_k |{\hat f(k)}+f(k)|^2\big).
        \label{Eq:MSE}
\end{equation}

\subsection{Noise-free reconstruction}\label{SubSec:NoNoise_Sim}

\begin{figure}[t]
\centering
      \includegraphics[width=0.45\linewidth]{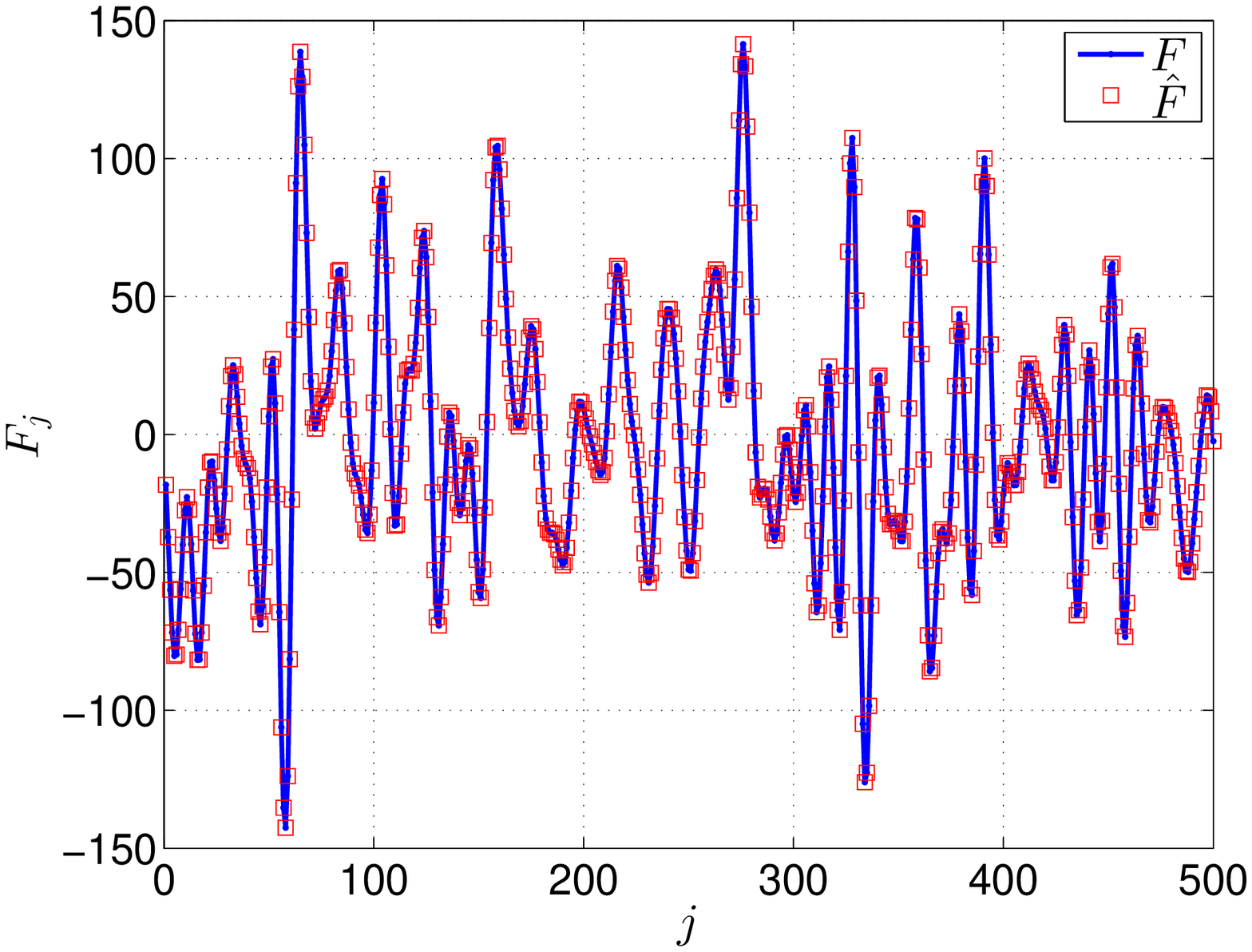}\ \
         \includegraphics[width=0.43\linewidth]{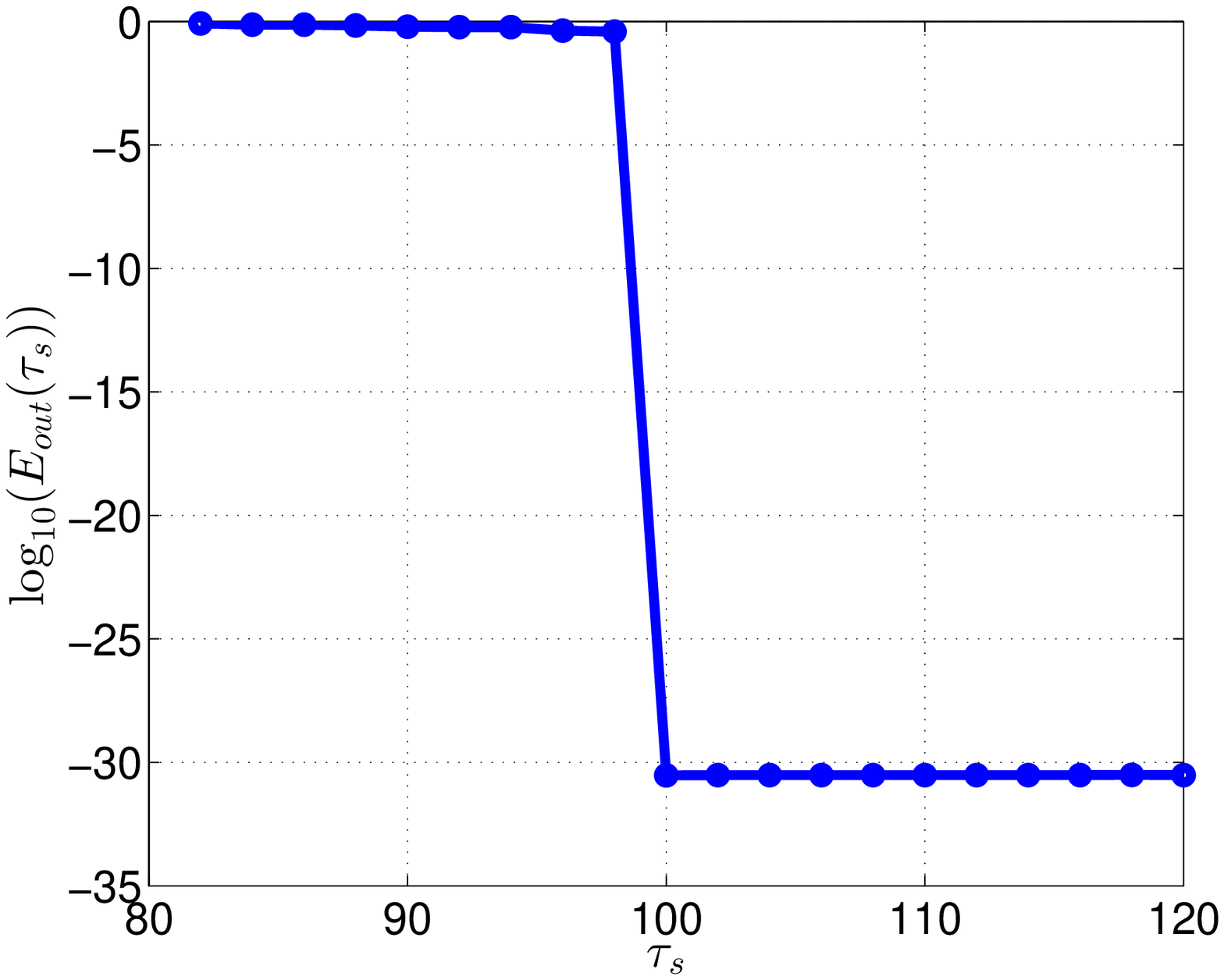}  
       
     \caption{Reconstruction of a signal with compact support \(\tau=100\), from \(N=500\) noise-free measurements. (left) The reconstructed versus true signal;  (right) the residual energy outside the compact support, on a log-scale, as a function of the unknown compact support. Our algorithm correctly estimates $\hat\tau=100$, and perfectly recovers the unknown sign pattern, hence the  MSE is essentially zero.}            
                \label{Fig:NoNoise_Rec}
\end{figure}

We demonstrate that with noise-free measurements and sufficient over-sampling, our sign retrieval algorithm perfectly recovers the underlying signal, essentially with machine-precision error. To this end, we generated a complex-valued random signal of length \(N=500\) with compact support parameter \(\tau=100\), whose Fourier transform is real-valued. We applied our sign retrieval algorithm including a scan over the unknown compact support, as described in Section \ref{Section:SignAlgo}. The left panel of Fig.  \ref{Fig:NoNoise_Rec}
shows the reconstructed signal versus the true one in the frequency domain. The right panel shows, on a logarithmic scale, the residual energy outside the assumed compact support, as described in Section \ref{SubSec:tau_estimation}. As seen from this figure and in accordance with our theoretical analysis, the residual error sharply drops precisely at the correct compact support parameter $\tau_s=100$ and remains at essentially zero value (up to machine error) for all values $\tau_s \ge \tau$. At this estimated compact support, our algorithm recovers correctly the exact sign pattern, leading to a zero MSE. 

\begin{figure}[t]
\centering
      \includegraphics[width=0.43\linewidth]{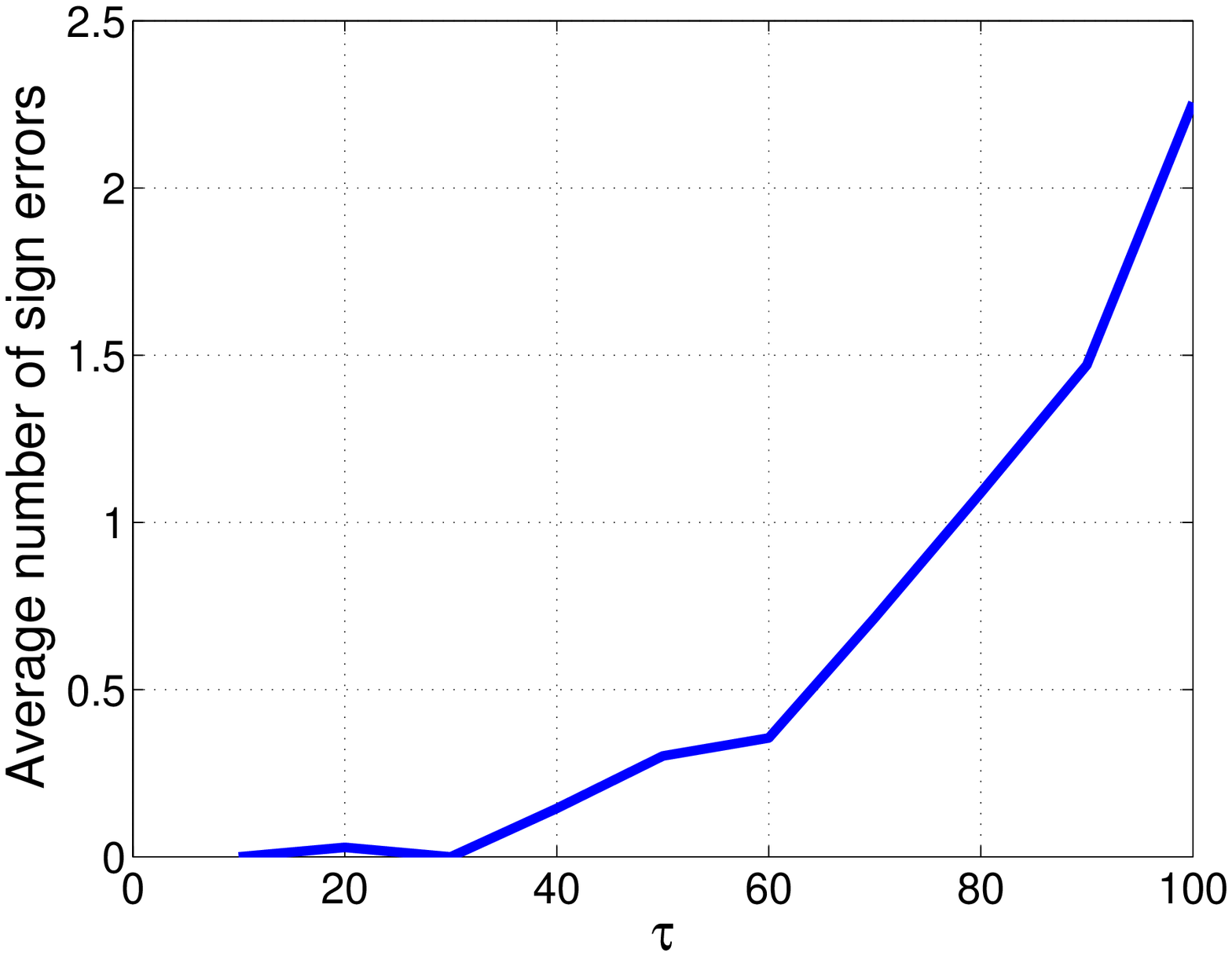}\ \
         \includegraphics[width=0.43\linewidth]{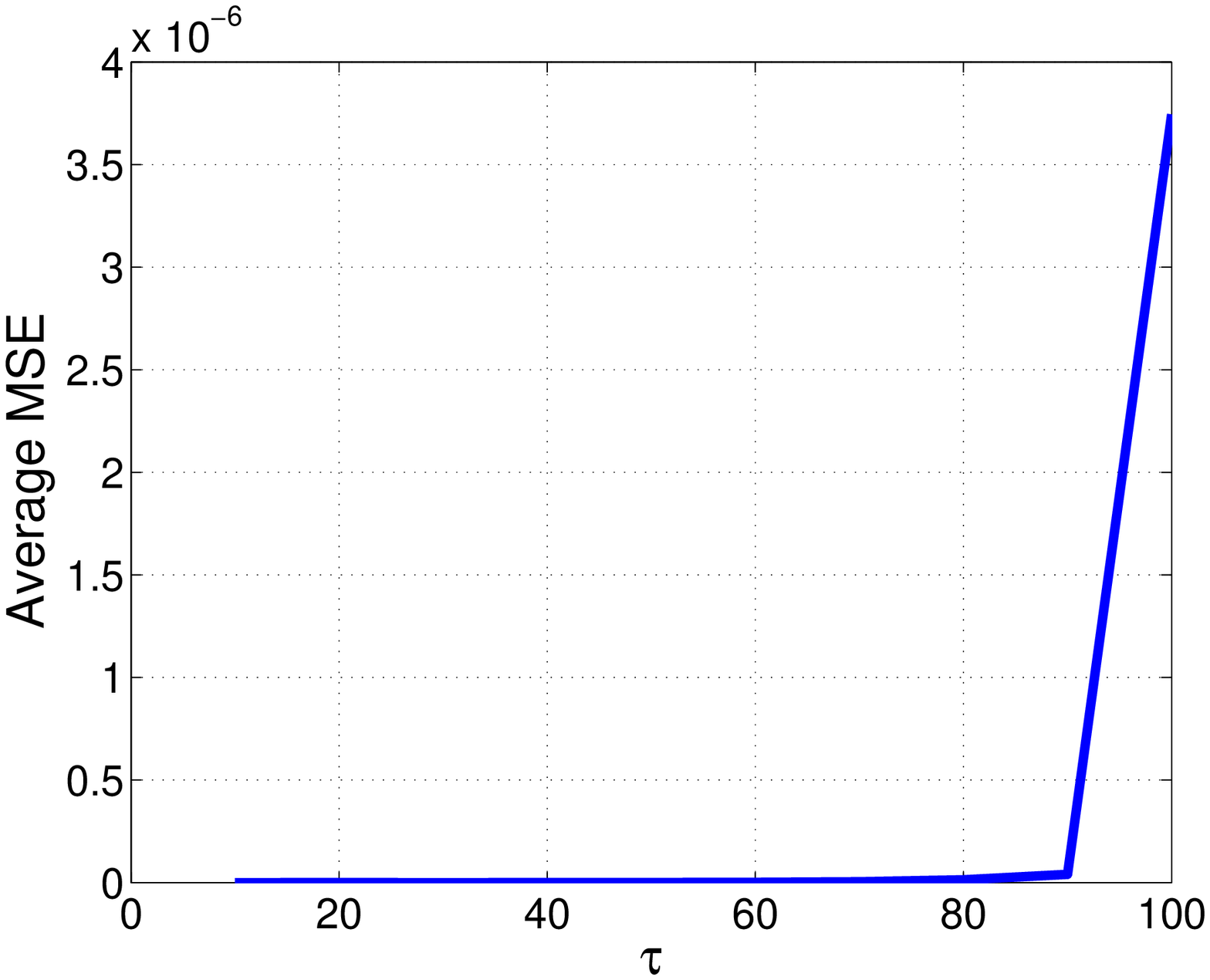}  
       
     \caption{Reconstruction of a random signal as a function of the compact support parameter $\tau$ with $N=500$. (Left) Average number of sign errors over 100 realizations. (Right) Average MSE.}            
                \label{Fig:Rec_Pr}
\end{figure}

In the above example, our algorithm found a correct over-segmentation with not too many segments and thus resulted in perfect reconstruction. Even with noise-free measurements, this is not always the case. In the left panel of Fig. \ref{Fig:Rec_Pr}, we present the number of sign errors, averaged over 100 random realizations as a function of the compact support parameter \(\tau\), with a fixed signal length $N=500$, whereas the right panel shows the corresponding averaged MSE. In accordance with Remark \ref{Remark:Large_N}, for small values of $\tau$ (e.g., a high oversampling rate), our segmentation scheme almost always obtains a perfect recovery. In contrast, as \(\tau\) increases, the exact segmentation does not yield a sufficient number of equations, and the heuristic segmentation may make small errors. These, however, typically occur at small values $|F_j|$, which as seen in the right panel lead to small reconstruction errors.

\subsection{Reconstruction in the presence of noise}\label{SubSec:Noisy_Sim}
\begin{figure}[t]
     \centering
      \includegraphics[width=0.45\linewidth]{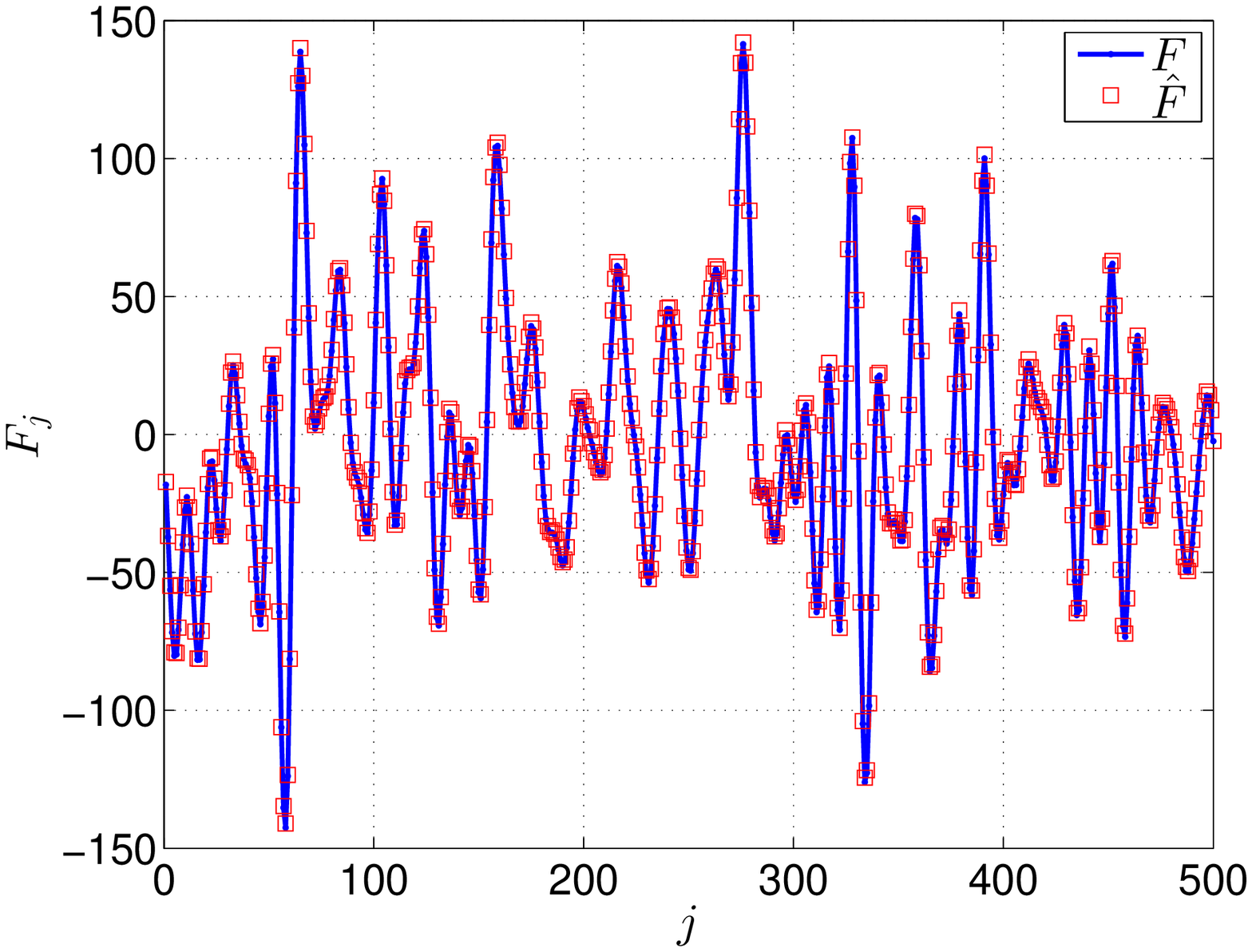}\ \ \
         \includegraphics[width=0.425\linewidth]{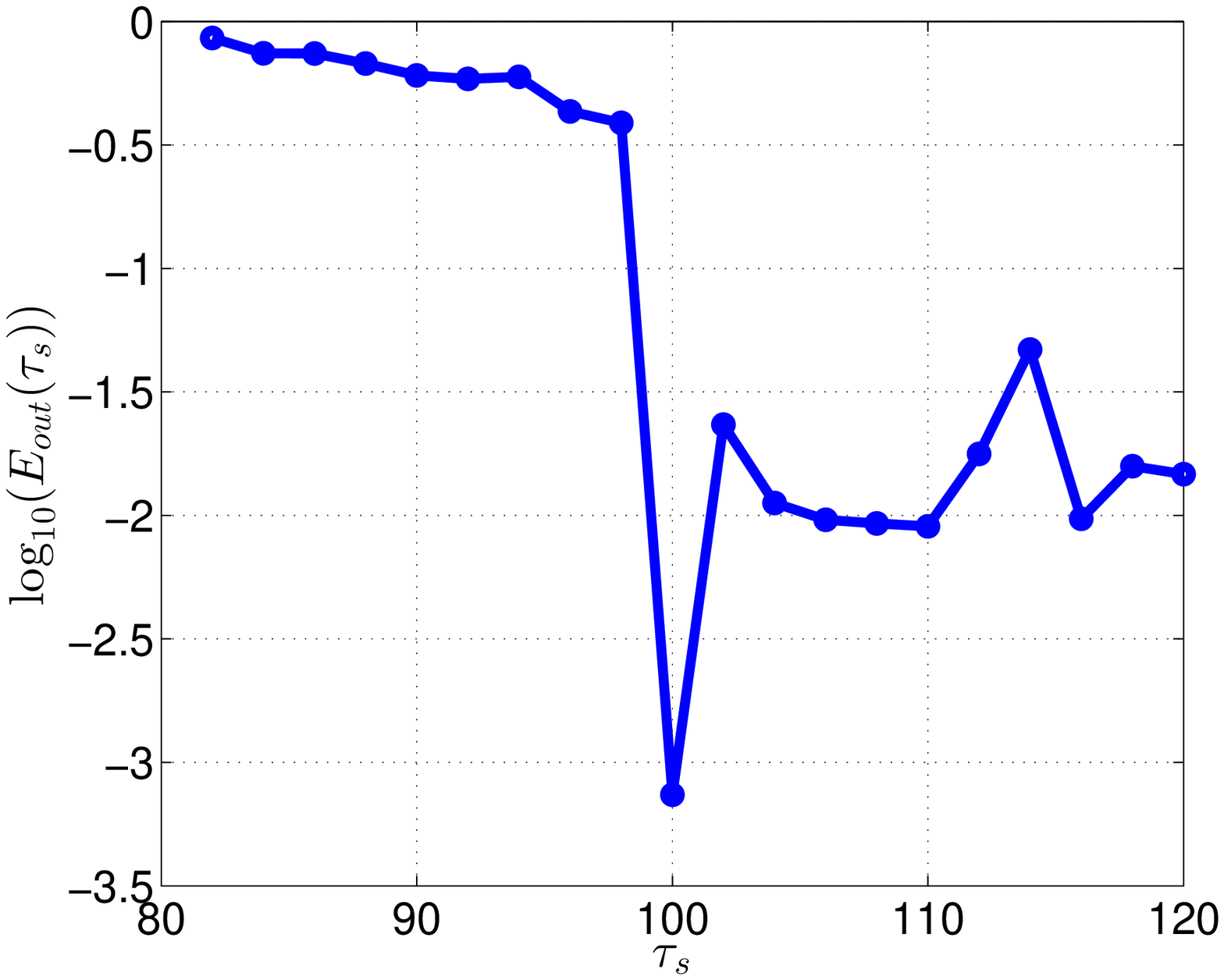}
     \caption{Reconstruction of the same signal as in Fig. \ref{Fig:NoNoise_Rec} in the presence of noise at level $\sigma=0.03$. (Left) Signal reconstruction. (Right) Residual error vs. assumed compact support \(\tau_s\).}
        \label{Fig:Noisy_Rec}
\end{figure}

Next, we demonstrate the robustness of our algorithm to noise. First, we show the resulting reconstructions from noisy measurements with  $\sigma=0.03$ of the same signal as in Section \ref{SubSec:NoNoise_Sim}, normalized to  $\|{\bf F}\|=1$. As shown in Fig. \ref{Fig:Noisy_Rec}(left) the reconstruction is in very good agreement with the true signal. The scan over $E_{out}$ for different values of the compact support depicted in Fig. \ref{Fig:Noisy_Rec}(right), shows that in contrast to the noise-free case, $E_{out}$ increases as the compact support is scanned above its true value. This occurs because increasing the compact support is equivalent to decreasing the number of linear equations which in turn decreases the stability to noise of the minimization problem described in Section \ref{SubSec:SR_Alg}.


%%%%%%%%%%%%%%%%%%%%%%%%%%
%

\begin{figure}[t]
        \centering
        \includegraphics[width=0.45\linewidth]{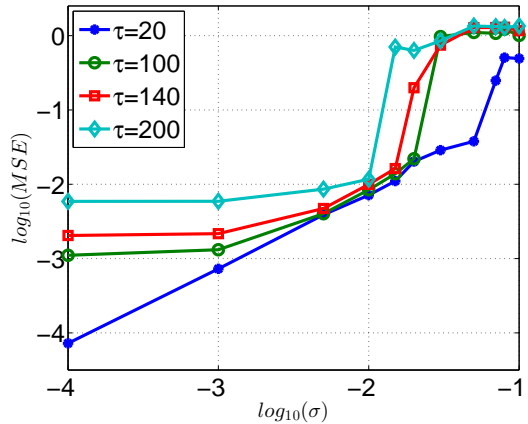}
        \caption{Reconstruction error versus noise level (log-log scale).}
        \label{Fig:MC_4taus}
\end{figure}

Next, we study the effect of noise on our reconstructions  by Monte Carlo simulations for different compact support lengths, $\tau = 20, 100, 140,200$, while keeping the oversampling ratio constant at  $N=5\tau$. For each $\tau$, we computed the average MSE from the reconstructions of $10$ randomly generated, complex-valued signals, each with $100$ noise realizations. The results, presented in  Fig. \ref{Fig:MC_4taus}, show that the reconstructions are stable up to $\sigma \simeq 0.01$ even for $\tau=200$. Note that for $\tau=20$, the segmentation scheme makes almost no errors, and thus the MSE increases linearly with $\sigma$ (on a log-log scale). For higher values of $\tau$, the probability for small segmentation errors increases, leading to nearly constant MSE for small noise levels. Possible routes to improve the noise stability are discussed in Section \ref{Sec:discussion}. 

\subsection{Reconstruction of a continuous function}
\label{SubSec:Compre_Thakur}

We now study the applicability of our algorithm, which relies on a discrete formulation, to the continuous sign problem. Following Thakur \cite{Thakur}, we sample at $N=200$ equispaced points the absolute $|F_c(\omega_j)|$ of a real-valued continuous function $F_c(\omega)$ whose continuous Fourier transform $f_c(t)$, has a compact support. Note that this scenario is different from our previous investigations, since the discrete Fourier transform of $F_c(\omega_j)$ is \textit{not} compactly supported. It is thus interesting to see if our algorithm can still succeed in recovering the sign pattern.

\begin{figure}[t]
     \centering
      \includegraphics[width=0.45\linewidth]{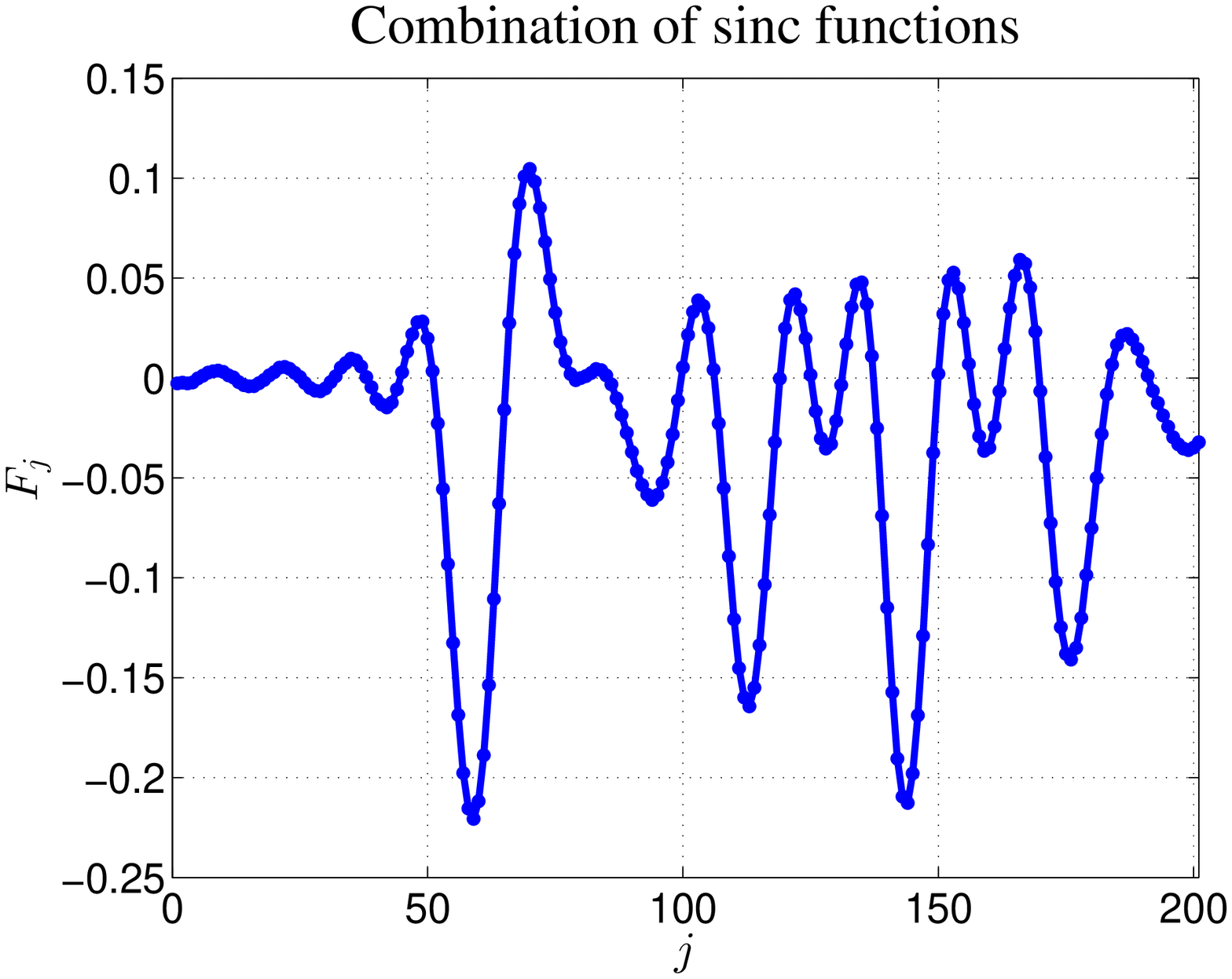}
       \includegraphics[width=0.45\linewidth]{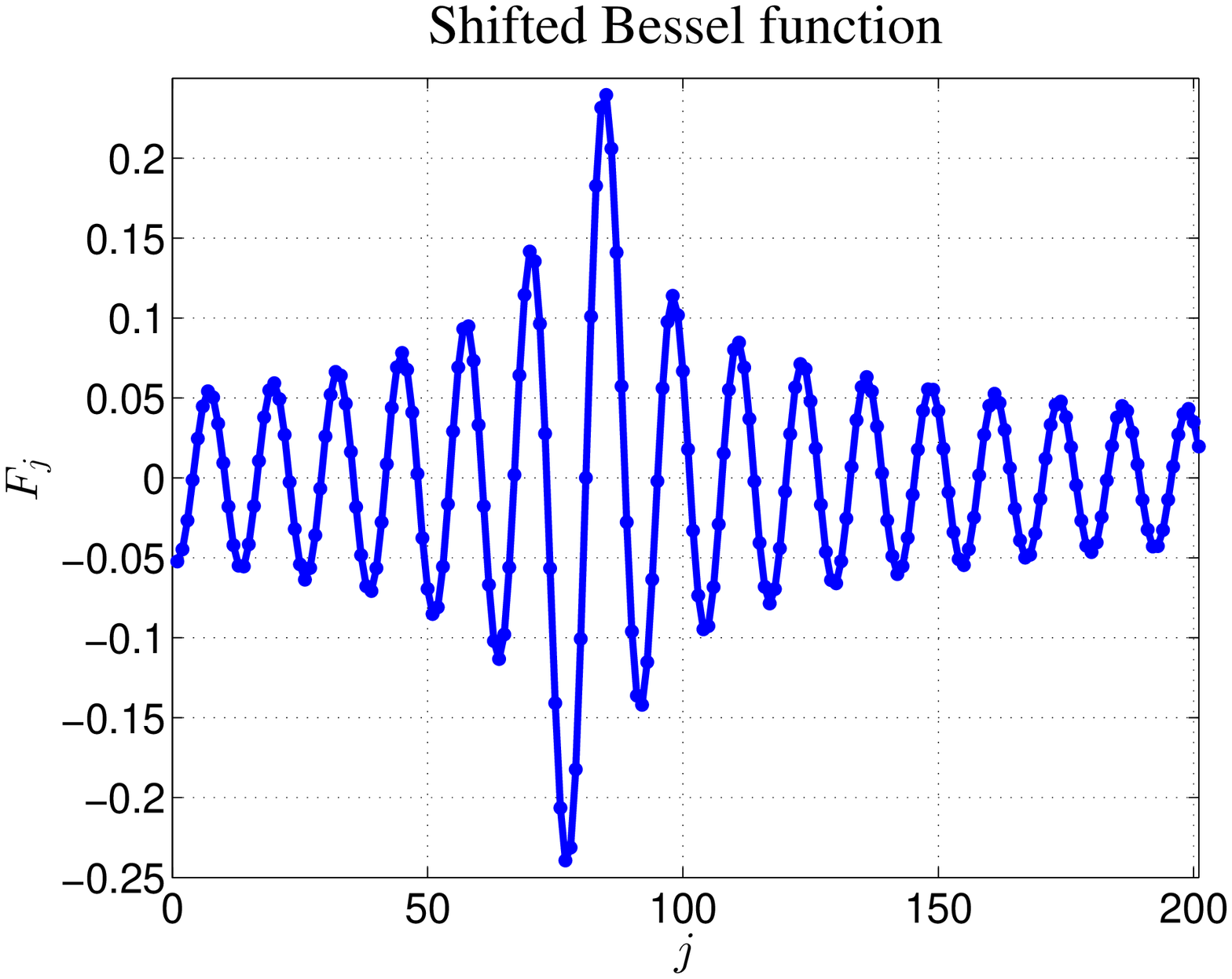}
     \caption{Two different real-valued continuous signals with band-limited Fourier Transform. }
     \label{Fig:Compare_Thakur_Funcs}
\end{figure}

\begin{figure}[t]
     \centering
     \includegraphics[width=0.45\linewidth]{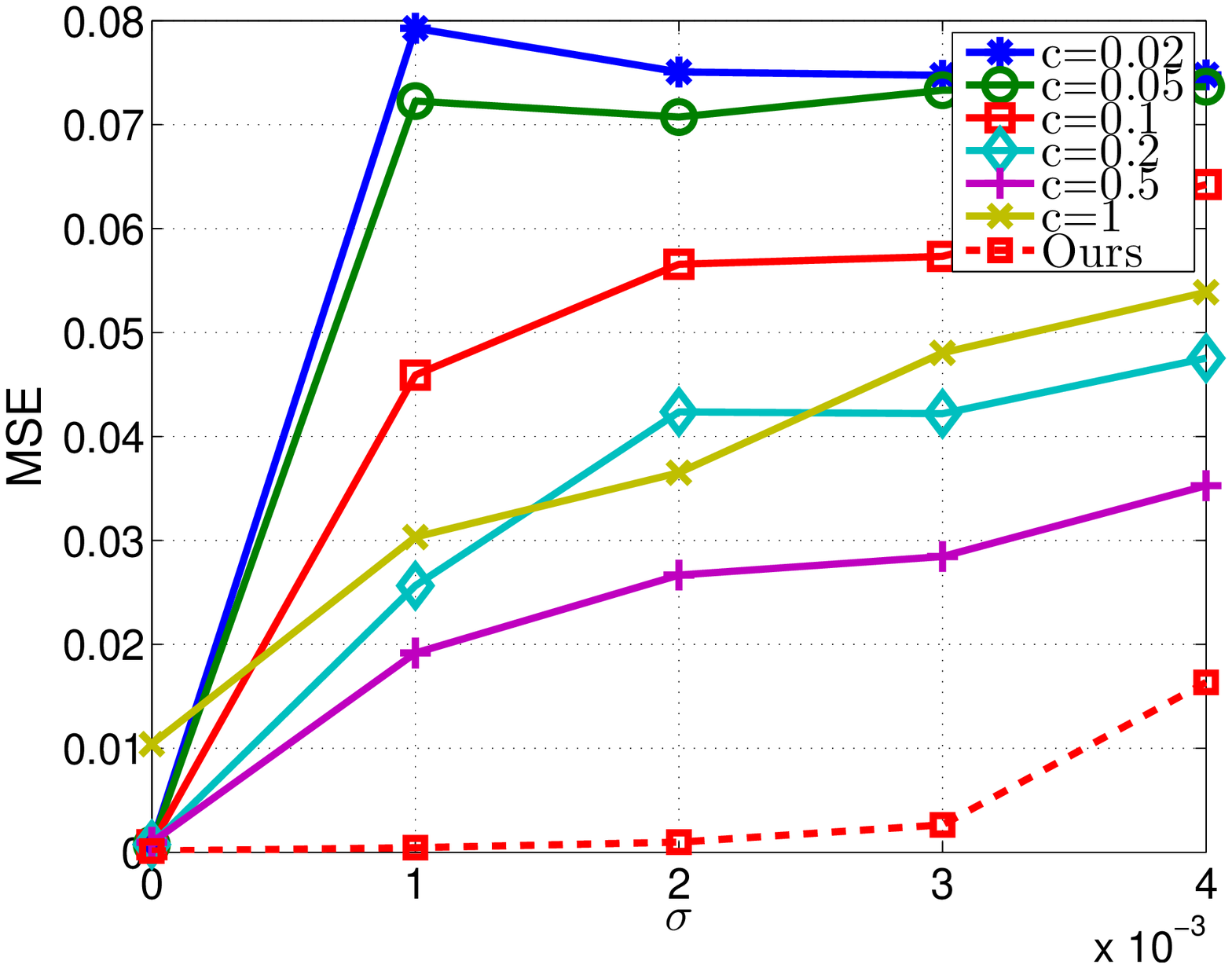}
      \includegraphics[width=0.45\linewidth]{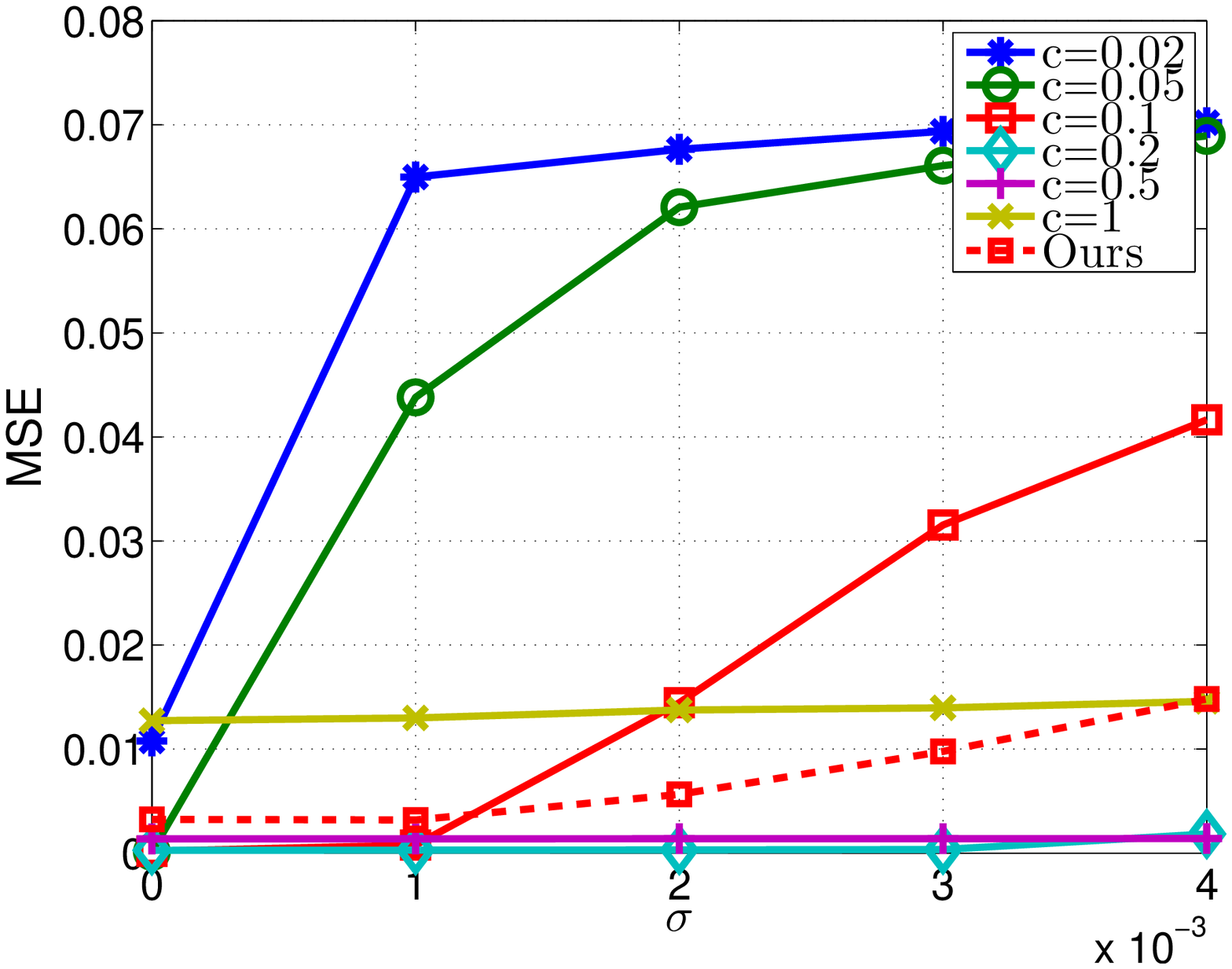}
     \caption{MSE comparison of Thakur's method with different values of \(c\) and ours as a function of noise level. (left) combination of sinc functions; (right) Shifted Bessel function. }
     \label{Fig:Compare_Thakur_MC}
\end{figure}

We compare our algorithm to \cite{Thakur} via Monte Carlo simulations at several noise levels, $\sigma = (0,1,2,3,4)\cdot10^{-3}$. We consider two signals, a shifted Bessel function as in \cite{Thakur} and depicted in Fig. \ref{Fig:Compare_Thakur_Funcs}(right), and a linear combination of 10 randomly shifted, randomly weighted sinc functions, depicted in Fig. \ref{Fig:Compare_Thakur_Funcs}(left). 
Note that Thakur's method depends on a parameter $c$ that needs to be set. Following the recommendation in \cite{Thakur}, we considered the following six values for reconstruction, $c=0.02,0.05,0.1,0.2,0.5,1$. 
Fig. \ref{Fig:Compare_Thakur_MC}(left) compares the two methods for the randomly shift sinc functions. In this case,  Thakur's method exhibits poor sensitivity to noise and its MSE rapidly increases with noise level for all considered values of $c$. In contrast, our method achieved a  significantly lower error. For the shifted Bessel function, presented in Fig. \ref{Fig:Compare_Thakur_MC}(right), Thakur's method achieved a lower MSE than ours for $c=0.2,0.5$. However, we note that the choice of the optimal $c$ value is not a-priori known, and may depend on the signal to be reconstructed. Also, as we verified in additional simulations, our method achieved a lower MSE than Thakur's, if the sampling rate is increased.

\subsection{VPR with 3 measurements}\label{SubSec:VPR_3M_Sim}

\begin{figure}[t]         
        \centering
        \includegraphics[width=0.45\linewidth]{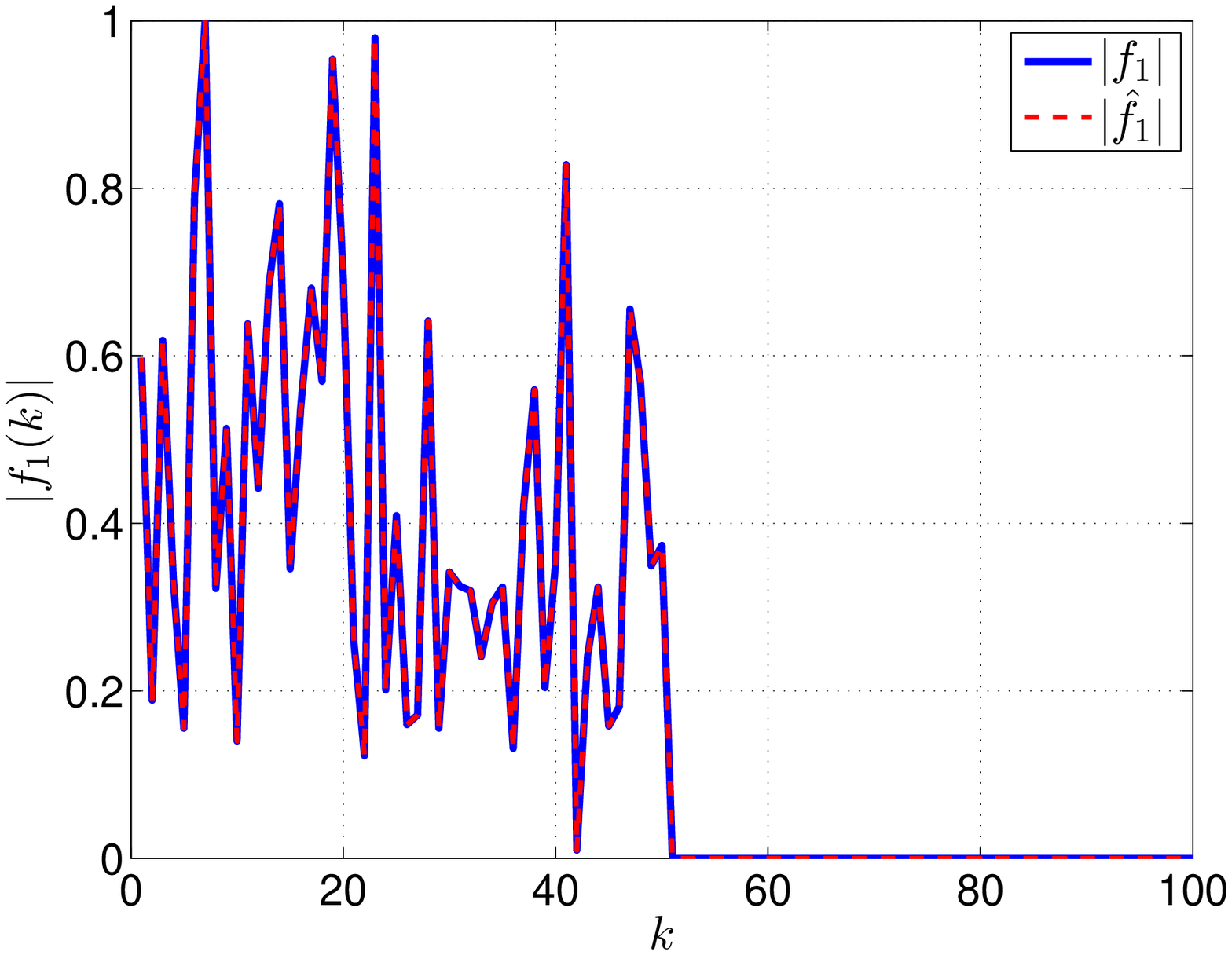}
        \includegraphics[width=0.425\linewidth]{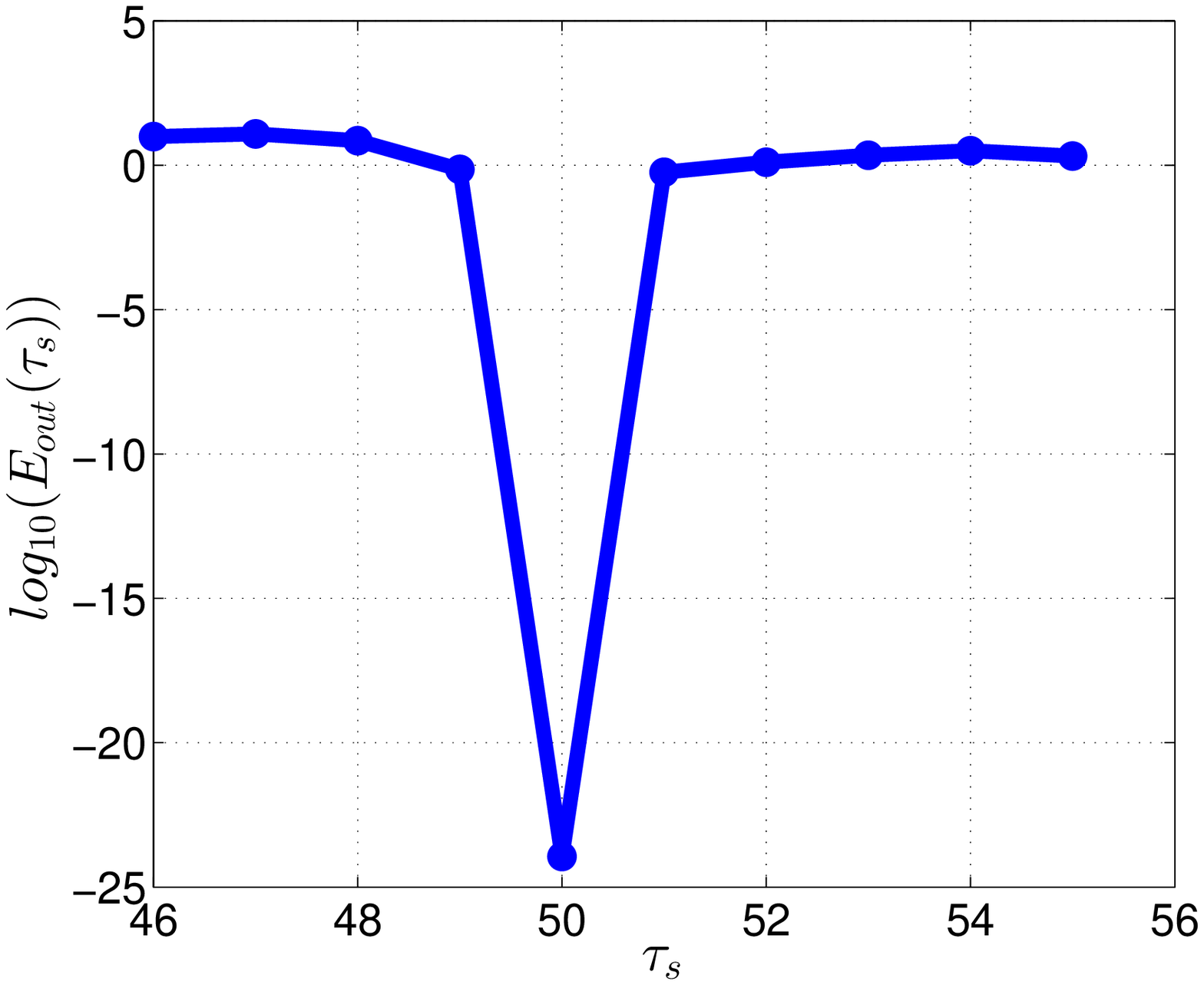}
        \caption{VPR with 3 measurements. (Left) Normalized magnitude of the reconstructed and true signal. (Right) Scan over averaged energy outside of the compact support. }
        \label{Fig:VPR_3M}
\end{figure}

Next, we demonstrate the use of sign retrieval with VPR to reconstruct two unknown complex-valued random signals, ${\bf f}_1$ and ${\bf f}_2$ both with compact support parameter $\tau=50$, from 3 measurements $|{\bf F}_1+{\bf F}_2|^2$, $|{\bf F}_1|^2$ and $|{\bf F}_2|^2$, as described in Section \ref{SubSec:VPR_3M}.  
Fig. \ref{Fig:VPR_3M}(left) shows the result for a noise-free  reconstruction. For simplicity, we present the magnitudes of  the complex-valued true and reconstructed signal, ${\bf f}_1$ (the results for ${\bf f}_2$ are similar). The reconstruction is perfect to within machine error. In this reconstruction we do not assume that the compact support size is known. Instead, we scan over its possible values, similarly to Section \ref{SubSec:tau_estimation} and \cite{OurIEEE}. A plot of the average energy outside of the compact support is presented in Fig. \ref{Fig:VPR_3M}(right). Note that in contrast to the sign problem here a clear minimum of $E_{out}$ is obtained at the correct value of the compact support. The reason for this behavior is that for the phase problem, as opposed to the sign problem, shifted solutions are allowed as detailed in \cite{OurIEEE}.

\subsection{Phase retrieval of two separated objects}\label{SubSec:VPR_SepObjs_Sim}

\begin{figure}[t]         
        \centering
        \includegraphics[width=0.45\linewidth]{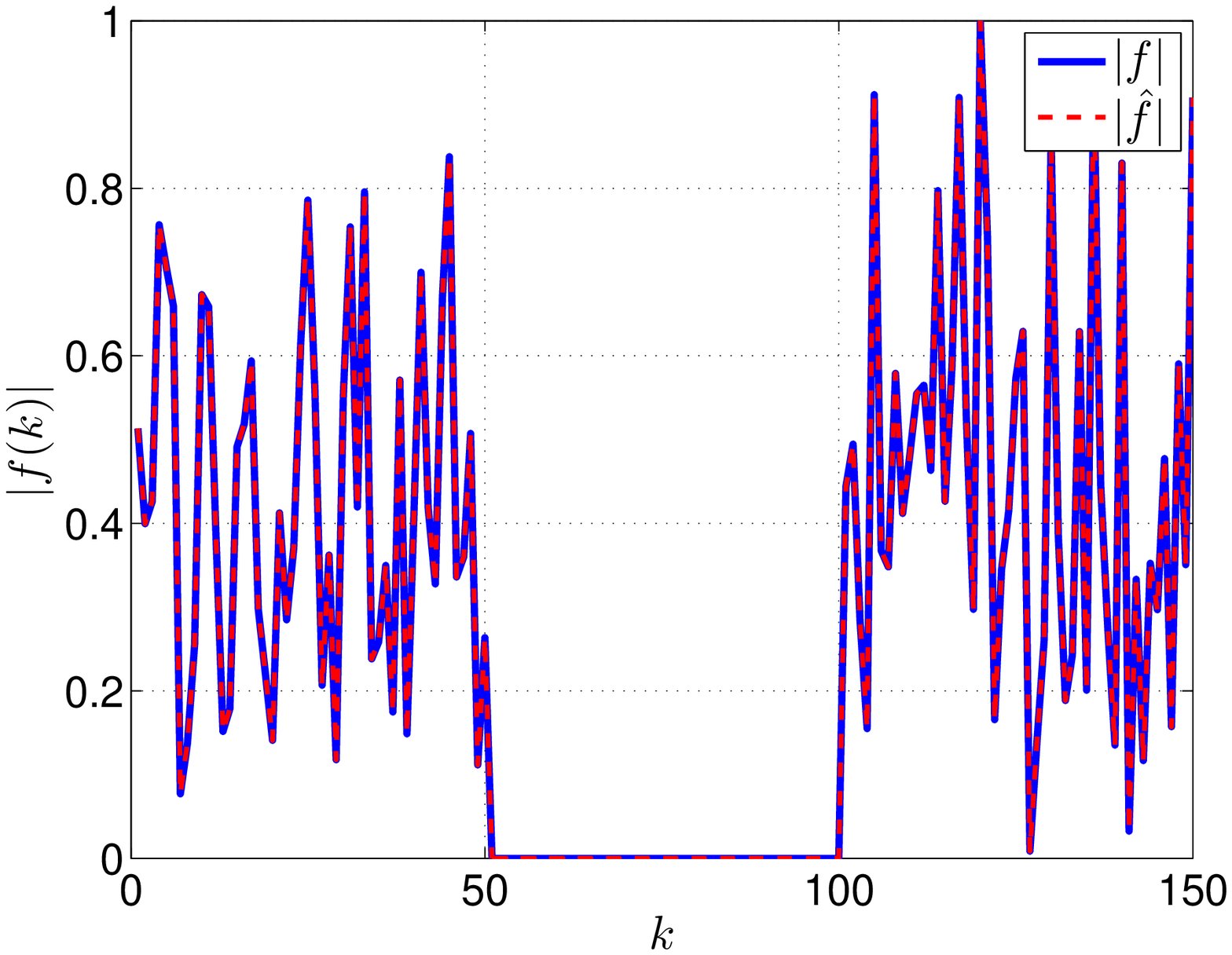}
        \includegraphics[width=0.425\linewidth]{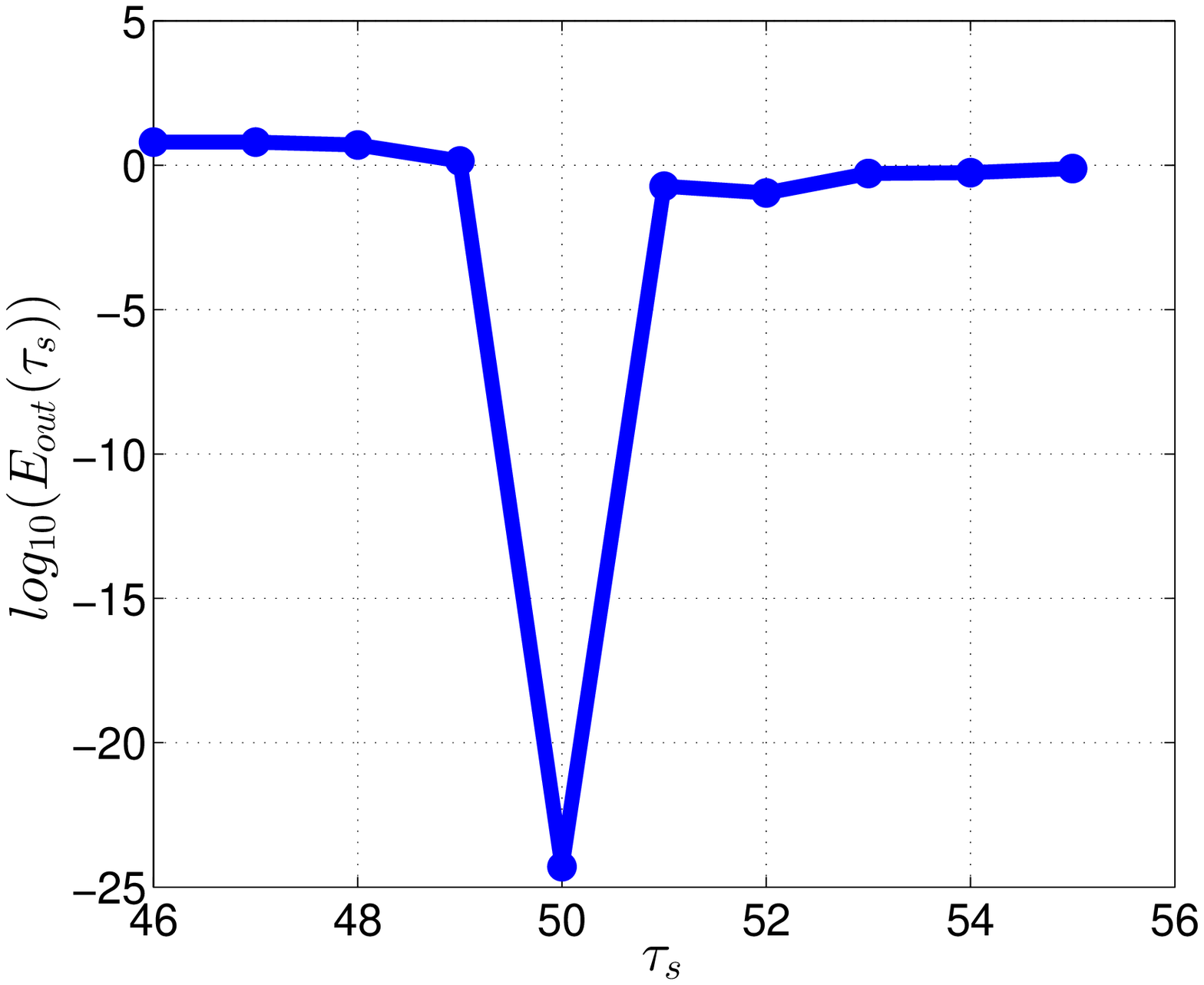}
        \caption{Phase retrieval with separated objects. (Left) Normalized magnitude of the reconstructed and true signal. (Right) Scan over averaged energy outside of the compact support.  }
        \label{Fig:SepObjs}
\end{figure}

Here we demonstrate 1D phase retrieval of a complex-valued vector consisting of two well separated signals, from a single measurement. This is achieved by combining sign retrieval with VPR as described in Section \ref{SubSec:VPR_SepObjs}. To this end, we generated a signal \({\bf f}\) of length \(N=500\) and compact support length  $151$, which consists of two randomly drawn, complex-valued vectors of length $50$ separated by $51$ zero valued entries. In the absence of noise, the signal is perfectly reconstructed (to within machine error), as demonstrated in Fig \ref{Fig:SepObjs}(left). As in section \ref{SubSec:VPR_3M_Sim}, we did not assume here a known compact support but rather estimated it as a part of the algorithm. The scan over $\tau_s$ for different compact support values is presented in Fig. \ref{Fig:SepObjs}(right). As mentioned in Section \ref{SubSec:VPR_3M_Sim}, also here a clear minimum of $E_{out}$ is attained at the true compact support value. 

%Note that although sign retrieval is used here as a step in the algorithm, $\bf F$ is complex-valued. Hence, this is a demonstration of 1-D phase retrieval from a single measurement without any prior knowledge other than that the compact support consists of two well separated regions.

\subsection{Segmentation scheme}\label{SubSec:SegScheme}
\begin{figure}[h]         
        \centering
        \includegraphics[width=0.8\linewidth]{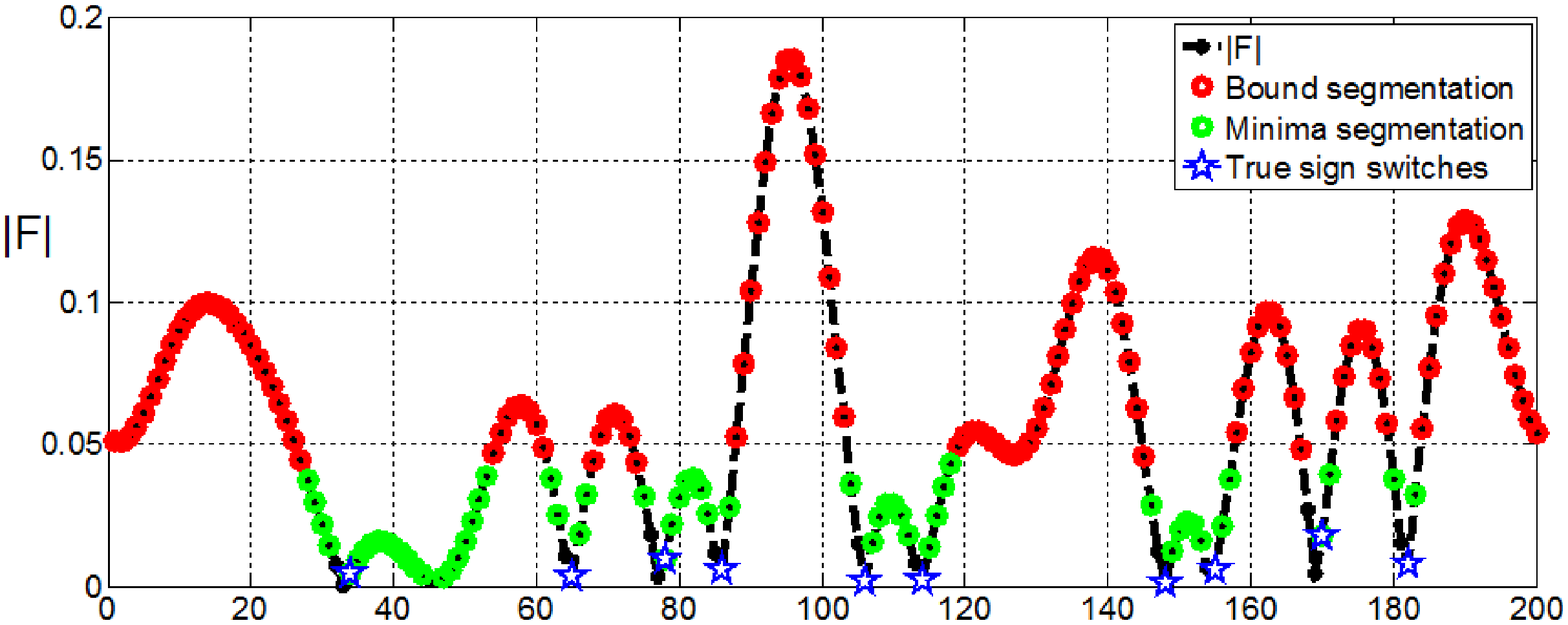}
        \caption{Illustration of the segmentation scheme. The black line denotes the absolute value of the signal $|F(\omega_j)|$. The red circles mark the regions guaranteed to have a common sign according to Lemma \ref{lem:CrossTh}. The green circles mark the additional regions that have common sign according to our heuristic segmentation scheme. The blue stars denote the true sign switching indices.   }
        \label{Fig:Seg_Example}
\end{figure}

Fig. \ref{Fig:Seg_Example} illustrates our segmentation scheme of Section \ref{SubSec:Seg} in the noise-free case. The indices in which sign changes occur are denoted by blue stars. The red circle denote indices for which Eq. \eqref{Eq:BoundCond} holds. The green circles denote indices for which Eq. \eqref{Eq:BoundCond} does not hold, but are assigned to constant sign intervals according to the heuristic segmentation scheme of Table \ref{Tab:HeuSeg}. 

\section{Discussion}
\label{Sec:discussion}

In this paper we presented a theoretical study of the finite dimensional sign problem and developed a novel computationally efficient sign retrieval algorithm. We then showed its applicability to two phase retrieval problems of practical interest, VPR with 3 measurements and reconstruction of two well separated objects.

The sign problem with a compactly supported underlying signal, which can be formulated as the solution to Eq. (\ref{Eq:X_CS}), is a specific instance
of an under-determined system of linear equations with an integer solution (the sign pattern in our case). The question of uniqueness of the solution to such systems, was recently studied by \cite{mangasarian2011probability}. In their paper, \cite{mangasarian2011probability} also proposed a linear programming relaxation of this non-convex problem, and proved it to work with overwhelming high probability if the equations are \textit{random} and their number is more than half the number of unknowns. Unfortunately, this relaxation was unable to recover the correct sign pattern in our sign problem, probably due to the fact that our equations are far from being random. 

Due to Lemma \ref{lem:sign_changes}, which states that the sign pattern can have at most \(\tau\) sign changes, our sign problem can also be viewed as a specific instance of an under-determined set of linear equations whose solution is sparse (for example by a change of variables from ${\bf s}$ to its discrete derivatives). For this problem, \(L_1\)-type minimization schemes have been proposed and proven to recover the correct solution under various conditions. Unfortunately, a straightforward application of $L_1$ Lasso penalization was again unable to recover the correct sign pattern, even with noise-free measurements, unless the compact support parameter, \(\tau\) which also controls the sparsity, was extremely small compared to \(N\).

Theoretical understanding of the robustness to noise of our algorithm is still lacking and is an interesting topic for further investigation. Improving the noise robustness is also an interesting route of future research. In particular, coupling the powerful relaxation schemes discussed above and their underlying theoretical guarantees to our segmentation-based relaxation could potentially improve the overall robustness of our sign retrieval algorithm. 

%Computational difficulty of sign problem when \(N=2\tau\) in the regime where our segmentation method fails to provide sufficient number of equations. 
\paragraph{Acknowledgments.}
We thank Dan Oron and Nirit Dudovich for useful discussions. We thank Gaurav Thakur for interesting discussions and for providing us with the code for his algorithm. We also thank Uri Feige for pointing us to ref.  \cite{mangasarian2011probability}. The research of B.N. was supported in part by  grant 892/14 from the Israeli Science Foundation. O.R. acknowledges financial support from the James S. McDonnell Foundation.  

\bibliographystyle{plain}
\bibliography{SignBib}

\appendix

\section*{Appendix}
\renewcommand{\thesubsection}{\Alph{subsection}}
\numberwithin{equation}{subsection}

\subsection{Proofs}\label{Appendix_Proofs}

\begin{proof}[\bf{Proof of Theorem \ref{Theorem:SignUniq}}]

By assumption A1, the compact support of the signal \({\bf f}\) is of the form $[-\tau/2,\ldots,\tau/2]\mod N$. It is convenient to consider its circularly {\em shifted} signal, ${\bf f}_s$, defined as $f_s(k)=f(k+\tau/2\mod N)$. By definition, ${\bf f}_s$ has compact support at the indices $0,1,\ldots,\tau$, and thus its $\mathcal{DFT}$ is by definition
$F_s(\omega_j) = \sum_{k=0}^\tau f_s(k)e^{\i\omega_j k}$. 
The relation between $\bf F_s$ and $\bf F$ is 
\[
F(\omega_j)=e^{-\i \omega_j\tau/2}F_s(\omega_j) = e^{\i\omega_j (N- \frac{\tau}{2})} F_s(\omega_j)
\]
where the last equality follows since $e^{\i \omega_j N}=1$ for all $j$. 

Next, we make the change of variables  $z=e^{\i\omega}$, known as the z-transform. With some abuse of notation we denote the resulting polynomial as $F(z)$, defined for all $z\in\mathbb{C}$,
\begin{equation} \label{Eq:E_Z}
F(z)=z^{N-\frac{\tau}{2}}\sum_{k=0}^{\tau} f_s(k)z^k.
\end{equation}
By its definition, $F(z)$ is analytic in the complex plane. Further, at the \(N\) sampling points, $z_j=e^{\i\omega_j}$, we have that $F(z)=F(\omega_j)$.

By the fundamental theorem of algebra, the polynomial 
$F(z)$ defined in Eq.\eqref{Eq:E_Z} can be decomposed as
\begin{equation}
F(z) = c z^{N-\frac{\tau}{2}} \prod_{r=1}^{\tau} (z-z_r) 
\label{eq:F_z}
\end{equation}
where $c \in \mathbb{C}$ is a normalization constant. Since the compact support of the shifted signal $f_s(k)$ is the set $\left\{0,1,...,\tau \right\}$ it follows that both $f_s(0)\neq 0$ and $f_s(\tau)\neq 0$. Hence, all $z_r \neq 0$. That is, the polynomial $F(z)$ has exactly $\tau$ non-zero roots. 

To guarantee a unique solution to the sign problem, these roots must be uniquely determined from the \(N\)
measurements $\{|F(\omega_j)|\}_{j=0}^{N-1}$. 
In the classical 1-D phase problem, where ${\bf F}\in\mathbb{C}^N$, these \(N\) measurements do not yield a unique solution. In our sign problem, in contrast, we have that ${\bf F} \in \mathbb{R}^N$, and as we now show, this leads to uniqueness, up to trivial ambiguities of the sign problem.  

To this end, consider the polynomial $F(z)^2$. By Eq. (\ref{eq:F_z}) it is given by
\begin{equation}
F(z)^2 = c^2 z^{2N-\tau} \prod_{r=1}^{\tau} (z-z_r)^2 
\end{equation}
Namely, $F(z)^2$ has the \textit{same} roots as $F(z)$, but each with its multiplicity doubled. Therefore, $F(z)^2$ has exactly $2\tau$ non-zero roots and is  thus uniquely determined by its values at any $2\tau+1$ distinct points.
Since $F(z_j)$ is real valued, at the $N$ observed points $F(z_j)^2=|F(z_j)|^2$. Thus, provided that  $N>2\tau$ these observations uniquely determine the \(2\tau\) roots, which in turn determine $F(z)$ up to a global \(\pm 1\) sign ambiguity.  
\end{proof}

 \begin{proof}[\bf{Proof of Lemma \ref{lem:sign_changes}}]
 The assumption that  ${\bf F}\in\mathbb{R}^N$ implies that 
 $f(k)=f^*(-k \mod N)$ for all $k=0,\ldots,N-1$. 
 This, in turn, implies that the  function  $F(\omega)=\sum_{k=-\tau/2}^{\tau/2}f(k)e^{\i\omega k}$, which extends \(F(\omega_j)=F_{j}\) to all \(\omega\in[0,2\pi]\), is real-valued for all $\omega$.  
 %the unit circle, $F(z)$ equals the discrete time Fourier transform of ${\bf f}$, that is,  $F(z=e^{\i \omega})=\sum_{k=0}^{\tau} f(k) e^{\i \omega k}$ with $\omega \in \mathbb{R}$.
 As in the proof of Theorem \ref{Theorem:SignUniq}, we consider the \(z\)-transform, \(z=e^{\i \omega}\), and the polynomial  \(F(z)\) of Eq. (\ref{Eq:E_Z}). Then, $F(z)$ is real-valued for all \(|z|=1\).  
 As in the proof of Theorem 1,  $F(z)$ has $\tau$ non-zero roots. Thus, in particular, $F(z)$ can vanish in at most $\tau$ points on the unit circle. Furthermore,  at the sampling points \(z_{j}=e^{\i\omega_j}\) we have  $F(z_j)=F_j$. Hence, in order for $F_j$ and \(F_{j+1}\) to have different signs, the continuous function $F(z)$ must have a zero crossing somewhere along the arc between $z_j$ and $z_{j+1}$. Therefore, the total number of zeros of $F(z)$ bounds the maximal number of sign changes in the vector \({\bf s}=sign({\bf F})\)  to be $\tau$. We note, that if at one of the sampling points $F(z_j)=0$ then $sign(F(z_j))$ is ill-defined. In this case we define $sign(F(z_j))=sign(F(z_{j-1}))$. If $F({z_0})=0$ we define $sign(F({z_0}))=1$.
 \end{proof}

\begin{proof}[\bf{Proof of Theorem \ref{Theorem:PiecePhase}}]
 By assumption A2, the signal \({\bf f}\) has a 
compact support of length $\tau+1$,
$$
CS=\left\{k:k\in [-\alpha,-\alpha+\tau] \mod N \right\}
$$ 
where $\alpha\in[0,N-1]$ is in general unknown. 

It will be convenient to work with the circularly shifted signal \({\bf f}_s\) given by $f_s(k) = f(k+\alpha \mod N)$, whose compact support consists of the indices \([0,1,\ldots,\tau]\). Its $\mathcal{DFT}$ is 
$F_{s}(\omega_j) = \sum_{k=0}^\tau f_s(k)e^{\i\omega_j k}$ and it is related to \({\bf F}\) via 
$F_s(\omega_j) = e^{\i\omega_j \alpha} F(\omega_j)$. 

As in the proof of Theorem \ref{Theorem:SignUniq}, consider the polynomial 
$F(z)=z^{N-\alpha}F_{s}(z)$
where $F_s(z)=\sum_{k=0}^{\tau} f_{s}(k)z^k$. In contrast to the sign problem, where \(F(z)\) was real-valued for all \(|z|=1\), in this case, where the phase is only assumed to be piecewise constant, \(F(z)\)
may in general be complex-valued. We thus write it as  
\begin{equation} \label{Eq:E_X}
F(z) = |F(z)|X(z).
\end{equation}
Assumption A2 can be expressed as $X(z_m\gamma^n)=a_m, \ m=1,\ldots,M,\ n=0,\ldots,N_m-1$ where $z_m=e^{\i \omega_{c(m)}}$ is the first point in each segment, $\gamma=e^{\i\Delta\omega}$ and $a_m \neq 0$ are unknown constants of unit modulus. Note that if at some interior point of a segment $|F(z_j)|=0$, its phase $X(z_j)$ is ill defined. In such a case we define it to be equal to the phase of the left-most point in that segment.  

To show that \(|{\bf F}|^2\) uniquely determines \(F(z)\) up to a global phase \(e^{\i \phi}\),  assume to the contrary that there exists another signal \({\bf g}\neq {\bf f}\) whose compact support is of the form $[-\alpha',-\alpha'+\tau]\mod N$, where possibly \(\alpha'\neq \alpha\). Its $\mathcal{DFT}$, ${\bf G }= \mathcal{DFT}\{{\bf g}\}$
satisfies that $|{\bf G}|=|{\bf F}|$ and it has a piecewise constant phase in the same \(M\) segments as \({\bf F}\) as defined in assumption A2. 

We denote the circular shift of ${\bf g}$ by \({\bf g}_s\), and the corresponding polynomials by $G(z)=z^{N-\alpha'}G_{s}(z)$  and $G_s(z)=\sum_{j=0}^\tau g_s(k)z^k$ respectively. Similarly to Eq. \eqref{Eq:E_X} we write  $G(z)=|G(z)|Y(z)$. Assumption A2 implies that $Y(z_m\gamma^n)=b_m$, where $b_m \neq 0$  are unknown constants of unit modulus.
 
Next, we define the following polynomial, where $\gamma=e^{\i \Delta \omega}$ 
\begin{eqnarray} \label{Eq:Pz}
P(z) &=& F(z) G(\gamma z)-G(z) F(\gamma z) \nonumber \\
     &=& z^{2N-\alpha-\alpha'}\left(\gamma^{N-\alpha'}F_s(z)G_s(\gamma z)-\gamma^{N-\alpha}F_s(\gamma z)G_s(z)\right).
\end{eqnarray}
Since both $F_s(z)$ and $G_s(z)$ are polynomials of degree $\tau$, the term inside the brackets in Eq. \eqref{Eq:Pz} is a polynomial of degree at most $2\tau$. Hence \(P(z)\) may have at most $2\tau$ non-zero roots. 

Now let us study the values \(P(z_{j})\). In each segment \(m\) of length \(N_m\geq 2\) we have that \(P(z)\) vanishes at all the points in this segment excluding its last one, $z_m \gamma^{N_m-1}$, since for any $n=0,\ldots,N_m-2$
\begin{align} \label{Eq:P_vanish}
\begin{split}
P(z_m \gamma^n) & = F(z_m \gamma^n) G(z_m \gamma^{n+1})-G(z_m \gamma^n) F(z_m \gamma^{n+1})\\
& = |F(z_m \gamma^n)||F(z_m \gamma^{n+1})|\Big(X(z_m \gamma^n) Y(z_m \gamma^{n+1})-Y(z_m \gamma^n) X(z_m \gamma^{n+1})\Big) \\
& = |F(z_m \gamma^n)||F(z_m \gamma^{n+1})|(a_m b_m-b_m a_m) = 0.
\end{split}
\end{align}
The total number of points where \(P\) vanishes in Eq.\eqref{Eq:P_vanish} is \(N-M.\) 
The condition that $N-M>2\tau$  implies that \(P(z)=0\) everywhere. 

Hence, specifically at the last point in each segment, upon division by the non-vanishing signal magnitudes\(\) 
\begin{equation}\label{Eq: XY1}
X(z_m \gamma^{N_m-1}) Y(z_m \gamma^{N_m})=Y(z_m \gamma^{N_m-1}) X(z_m \gamma^{N_m}), \  m=1, \ldots, M
\end{equation}
Next, using $z_m \gamma^{N_m}=z_{m+1}$ and the fact that $\bf X$ and $\bf Y$ are phase vectors, Eq. \eqref{Eq: XY1} can be written as
\begin{equation}\label{Eq: XY2}
X(z_m \gamma^{N_m-1})X^*(z_{m+1})=Y(z_m \gamma^{N_m-1})Y^*(z_{m+1}),\  m=1, \ldots, M.
\end{equation}
From Eq. \eqref{Eq: XY2} it follows that the phase difference between any pair of consecutive segments in $X(z_j)$ is equal to the phase difference between the corresponding segments in $Y(z_j)$. Hence, we have that 
\begin{equation}
X(z_j)=e^{\i\phi}Y(z_j), \ j=0,\ldots,N-1
\end{equation}
where $\phi$ is an arbitrary global phase. 
Therefore, $F(z)=e^{\i\phi}G(z)$ and since $F_j=F(z_j)$ we conclude that ${\bf F}$ is uniquely determined up to an arbitrary global phase.

  We note that if $F_j$ vanishes at the last and/or first indices in one or more segments,  then our proof can be modified as follows: Redefine the polynomial $P(z)$ of Eq.\eqref{Eq:Pz} as $P(z)=F(z) G(\gamma^2 z)-G(z) F(\gamma^2 z)$. For this modified $P(z)$, from Eq.\eqref{Eq:P_vanish} we have that $P(z)$ vanishes at a total of $N-2M$ points. Hence, if $N>2M+2\tau$ the proof can be completed in a straight forward manner. Furthermore if $|\bf F|$ vanishes at a number of consecutive points near the edges of some (or all) segments the proof can be modified in a similar manner.     
\end{proof}

\begin{proof}[\bf{Proof of Theorem \ref{Theorem:SignAlgo}}]
By definition, the true sign pattern \({\bf s}=sign({\bf F})\) is a solution of Eqs. (\ref{Eq:X_CS})-(\ref{Eq:X_ConstPhase}). We now show that all solutions are of the form \({\bf X}=c\,{\bf s}\) where $c\in\mathbb{C}$. 

%To this end, as in the proof of Theorem \ref{Theorem:PiecePhase}, consider the circularly shifted signal \({\bf f}_s\) given by $f_s(k) = f(k+\tau/2 \mod N)$, whose compact support consists of the indices \([0,1,\ldots,\tau]\). Its $\mathcal{DFT}$ is 
%$F_s(\omega_j) = \sum_{k=0}^\tau f_s(k)e^{\i\omega_j k}$ and it is related to \({\bf F}\) via 
%$F_s(\omega_j) = e^{\i\omega_j \tau/2} F(\omega_j)$. Its z-transform is given by the polynomial $F(z)=z^{N-\tau/2}F_{s}(z)=|F(z)|X(z)$ where $F_s(z)=\sum_{k=0}^{\tau} f_{s}(k)z^k$.

 Assume to the contrary that there exists another non-zero solution ${\bf Y}\neq c {\bf s}$, where \(Y_j\in\mathbb{C}\) is not necessarily of unit modulus, and yet  satisfies Eqs.\eqref{Eq:X_CS}-\eqref{Eq:X_ConstPhase}. Let ${\bf G}=|{\bf F}| {\bf Y}$ be the signal associated with this solution. By Eq. \eqref{Eq:X_CS},\   ${\bf g}=\mathcal {DFT}\{{\bf G}\}$ has the same (or smaller)\ compact support as that of $\bf f$. As in the proof of Theorem \ref{Theorem:PiecePhase}, we denote its shifted signal by ${\bf g}_s$ and its corresponding polynomial by $G(z)=z^{N-\tau/2}G_s(z)$. % Note that at the \(N\) sampling points, $G(z_j)=|F(z_j)|Y_j$. 

Next, as in Theorem \ref{Theorem:PiecePhase}, we define the polynomial $P(z)$ of Eq. (\ref{Eq:Pz})\ which has at most $2\tau$ non-zero roots. Since by Eq. (\ref{Eq:X_ConstPhase}) the vector\ \({\bf Y}\) is piecewise constant in the same segments as $\bf X$, it follows from Eq.\eqref{Eq:P_vanish} that $P(z)$ vanishes at a total of $N-M$ points. Hence, the condition  $N-M>2\tau$  implies that \(P(z)=0\) everywhere. In particular, \(P(z_{j})=0\) for all $j=0,\ldots,N-1.$ Specifically for the last point at each segment, upon division by \(|F(z_{m}\gamma^{N_m-1})F(z_{m+1})|\), which is non-zero by our assumption,  we have
\begin{equation}\label{Eq:XY_Ratio}
\frac{X(z_m \gamma^{N_m-1})}{X(z_{m+1})} = \frac{Y(z_m \gamma^{N_m-1})}{Y(z_{m+1})}, \  m=1, \ldots, M-1
\end{equation}
Eq.\eqref{Eq:XY_Ratio} implies that the proportionality constants between each pair of consecutive segments in $\bf X$ and $\bf Y$ are equal. Hence, $Y_j=cX_j$ for all $j=0, \ldots,N-1$ for some constant $c \in \mathbb{C}$.
\end{proof}

\begin{proof}[\bf{Proof of Lemma \ref{lem:CrossTh}}]
        As in the proof of Lemma \ref{lem:sign_changes}, let $F(\omega)$ be the extension of $F_j=F(\omega_j)$ to all $\omega \in [0,2\pi]$. Since ${\bf F} \in \mathbb{R}^N$, $\bf f$ is symmetric-conjugate and it is also compactly supported. Hence, $F(\omega)$ can be expressed as        
        \begin{eqnarray}
        F(\omega) &=& \sum_{k=-\tau/2}^{\tau/2} f(k)e^{-\i \omega k} = f_0+\sum_{k=1}^{\tau/2} \big( f(k) e^{-\i \omega k}+f^\ast(k) e^{\i \omega k}\big) \notag \\ &=& f_0+ \sum_{k=1}^{\tau/2} 2 \mathcal R(f(k) e^{-\i \omega k}) 
        =  f_0+2\sum_{k=1}^{\tau/2} |f(k)|\cos(\omega k+\theta_k),
        \end{eqnarray}
        where $f_0=f(k=0)$.
        Since $\omega_{j}-\omega_{j-1}=\Delta \omega=\frac{2\pi}{N}$, 
\begin{equation}
|F(\omega_j)-F(\omega_{j-1})| \leq \frac{2\pi}N \max_\omega |\frac{d}{d\omega}F(\omega)|
\end{equation}        
Combining the last two equations gives        
       \begin{eqnarray}
        |F(\omega_j)-F(\omega_{j-1})| &\leq& \frac{4\pi}{N} \sum_{k=1}^{\tau/2} |f(k)| k      \notag
        \end{eqnarray}
Finally, by the Cauchy-Shwartz inequality, 
\begin{eqnarray}   
      |F(\omega_j)-F(\omega_{j-1})|  &\le& \frac{4\pi}{N}  \sqrt{\sum_{k=1}^{\tau/2} |f(k)|^2} \ \sqrt{\sum_{k=1}^{\tau/2}k^2}
        \label{Eq:Drvt_Bound}
        \end{eqnarray} 
        Next, using Parseval's theorem we have
        \begin{equation}\label{Eq:BoundSum_f}
                \sum_{k=1}^{\tau/2}|f(k)|^2=\frac{1}{2}\sum_{k=-\tau/2}^{\tau/2}|f(k)|^2-\frac{1}{2}|f_0|^2=\frac{\|{\bf F}\|^2}{2N}-\frac{1}{2}|f_0|^2\leq\frac{\|{\bf F}\|^2}{2N}
    \end{equation}
    Combining Eqs. \eqref{Eq:Drvt_Bound}-\eqref{Eq:BoundSum_f} yields Eq. \eqref{Eq:LemmaBound} which concludes the proof.
\end{proof}

\subsection{Segmentation in the presence of noise}
Here we describe a simple modification to Eq. \eqref{Eq:BoundCond} to account for noise. As described in Section \ref{SubSec:SR_Alg}, our noise model is given by Eq. \eqref{Eq:NoiseModel}. After neglecting dark counts and shot noise, we are left with
\begin{equation}
|\tilde {\bf F}|=|{\bf F}+\tfrac{\sigma}{\sqrt{N}}{\bm \eta}^s|.
\end{equation}
Hence, assuming a high SNR we have,
\begin{equation}
|\tilde {\bf F}|=|{\bf F}|+\frac{\sigma}{\sqrt{N}}\mathcal{R}(e^{-\i {\bm \phi}}{\bm \eta}^{s}) +O\left(\frac{\sigma^2}{N}\right).
\end{equation}
Since we assumed ${\bm \eta}^{s} \sim \mathbb{c} \mathcal{N}(0,1),$ we may thus write ${\bm \eta}^{s}=\frac{{\bf a}+\i{\bf b}}{\sqrt{2}}$ where $\bf a$ and $\bf b$ are i.i.d.  $\mathcal{N}(0,1)$. Hence, $\mathcal{R}(e^{-\i {\bm \phi}}{\bm \eta}^{s})=\frac{1}{\sqrt{{2}}}({\bf a}\cos{\bm \phi}+{\bf b}\sin{\bm \phi})$ is zero-mean Gaussian and its variance is $\frac{1}{2}$. In the presence of noise, $|\tilde F_j|+|\tilde F_{j+1}|$ is thus approximately the correct value $| F_j|+| F_{j+1}|$, perturbed by a Gaussian with variance $\sigma^2/N$. To account for noise, we add this standard deviation $\frac{\sigma}{\sqrt{N}}$ to the threshold on the right hand side of Eq. \eqref{Eq:BoundCond}. Thus Eq. \eqref{Eq:BoundCond_Noisy} follows.

\end{document}